\documentclass{scrartcl}

\usepackage[pdftex,
pdfauthor={Boris Kramer and Karen Willcox},
pdftitle={Nonlinear Model Order Reduction via Lifting Transformations and Proper Orthogonal Decomposition},
pdfkeywords={Nonlinear model reduction, reduced-order modeling, lifting and variable transformations, PDE, reactive flows, combustion},
pdfproducer={Latex with hyperref},
pdfcreator={pdflatex}]{hyperref}

\usepackage[automark]{scrpage2}
\ihead[]{}
\ohead[]{}
\ifoot[]{}
\ofoot[]{} 
{\begin{trivlist}\item[]{\bfseries\sffamily Keywords: }\ }
	{\end{trivlist}}

\newcommand{\bit}{\begin{itemize}}
\newcommand{\eit}{\end{itemize}}
\newcommand{\Dam}{\mathcal{D}}
\newcommand{\spe}{\psi}
\newcommand{\bzero}{\mathbf{0}}

\newcommand{\x}{\mathbf{x}}

\renewcommand{\u}{\mathbf{u}}
\newcommand{\E}{\mathbf{E}}
\newcommand{\A}{\mathbf{A}}
\newcommand{\B}{\mathbf{B}}
\newcommand{\C}{\mathbf{C}}
\newcommand{\N}{\mathbf{N}}
\newcommand{\G}{\mathbf{G}}
\newcommand{\HH}{\mathbf{H}}
\newcommand{\V}{\mathbf{V}}
\newcommand{\X}{\mathbf{X}}
\newcommand{\I}{\mathbf{I}}

\newcommand{\aiaaappendix}[1]{\section*{Appendix: #1}}

\usepackage[utf8]{inputenc}
\usepackage{graphicx}
\usepackage[version=3]{mhchem}
\usepackage{longtable,tabularx}

\usepackage{amsmath}
\usepackage{amsfonts}
\usepackage{amssymb}
\usepackage{textcomp} 
\usepackage[gen]{eurosym}

\newtheorem{example}{Example}

 \usepackage{wrapfig}
 \usepackage{threeparttable}
 \usepackage{dcolumn}
  \newcolumntype{d}{D{.}{.}{-1}}
\usepackage{multicol}
\usepackage{nomencl}
\usepackage{etoolbox,ragged2e,siunitx,mathtools} 
 \usepackage{subfigure}
 \usepackage{fancyvrb}
  \fvset{fontsize=\footnotesize,xleftmargin=2em}
  \usepackage{fix-cm}
  \usepackage{setspace}
 \usepackage{lettrine}

 \usepackage{cleveref}

\renewcommand\nomgroup[1]{\def\nomtemp{\csname nomstart#1\endcsname}\nomtemp}
  
\renewcommand*{\nompreamble}{\begin{multicols}{2} \markboth{\nomname}{\nomname}}

\renewcommand{\nomname}{Nomenclature}
\renewcommand*\nompostamble{\end{multicols}}

\makenomenclature

\crefname{equation}{Eq.}{Equations}
\Crefname{equation}{Eq.}{Equations} 
\crefname{figure}{Fig.}{Figures}
\Crefname{figure}{Fig.}{Figures}
\crefname{table}{Table}{Tables}
\Crefname{table}{Table}{Tables}
\crefname{chapter}{Chapter}{Chapters}
\Crefname{chapter}{Chapter}{Chapters}
\crefname{appendix}{Appendix}{Appendices}
\Crefname{appendix}{Appendix}{Appendices}
\crefname{section}{Section}{Sections}
\Crefname{section}{Section}{Sections}
\crefname{algocf}{Algorithm}{Algorithms}
\Crefname{algocf}{Algorithm}{Algorithms}

\usepackage{pgfplotstable}      
\usepackage{booktabs,colortbl}  
\usepackage{array}
\usepackage{hhline}
\usepackage{float}
\usepackage{epstopdf}
\usepackage{wrapfig}
\usepackage[color=red,colorinlistoftodos]{todonotes}

\usepackage{tikz,pgfplots}
\usetikzlibrary{arrows,shapes,positioning,shadows,trees,patterns,intersections,matrix,fit,backgrounds}
\pgfplotsset{compat=newest} 
  \pgfplotsset{plot coordinates/math parser=false}
\newlength\figureheight 
    \newlength\figurewidth 
\usepackage{tikz-3dplot} 
\usepackage{tkz-euclide}
\usetkzobj{all}
\tikzset{
  basic box/.style={
    shape=rectangle, rounded corners, align=center,
    draw=#1, fill=#1!25},
}

\pgfplotstableset{
    every even row/.style={
    before row={\rowcolor[gray]{0.9}}},
    every even column/.style={
    column type={l}},
    every odd column/.style={
    column type={l}},
    every head row/.style={
    before row=\toprule,after row=\midrule},
    every last row/.style={
    after row=\bottomrule},
    col sep = &,
    row sep=\\,
    string type,
}

\definecolor{ChartColor1}{RGB}{236,231,242}
\definecolor{ChartColor2}{RGB}{166,189,219}
\definecolor{ChartColor3}{RGB}{43,140,190}
\newcommand{\add}[1]{{\color{black}~#1}}

\title{Nonlinear Model Order Reduction via Lifting Transformations and Proper Orthogonal Decomposition}

\author{
	Boris Kramer\thanks{Department of Aeronautics and Astronautics, Massachusetts Institute of Technology ({bokramer@mit.edu}).}
	\and Karen Willcox\thanks{Department of Aerospace Engineering and Engineering Mechanics and Institute for Computational Engineering and Sciences, University of Texas at Austin ({kwillcox@ices.utexas.edu}).}
}
\date{\today}

\begin{document}

\maketitle

\begin{abstract}
\noindent 
This paper presents a structure-exploiting nonlinear model reduction method for systems with general nonlinearities.
First, the nonlinear model is lifted to a model with more structure via variable transformations and the introduction of auxiliary variables. The lifted model is equivalent to the original model---it uses a change of variables, but introduces no approximations. When discretized, the lifted model yields a polynomial system of either ordinary differential equations or differential algebraic equations, depending on the problem and lifting transformation. Proper orthogonal decomposition (POD) is applied to the lifted models, yielding a reduced-order model for which all reduced-order operators can be pre-computed.
Thus, a key benefit of the approach is that there is no need for additional approximations of nonlinear terms, in contrast with existing nonlinear model reduction methods requiring sparse sampling or hyper-reduction.
Application of the lifting and POD model reduction to the FitzHugh-Nagumo benchmark problem and to a tubular reactor model with Arrhenius reaction terms shows that the approach is competitive in terms of reduced model accuracy with state-of-the-art model reduction via POD and discrete empirical interpolation, while having the added benefits of opening new pathways for rigorous analysis and input-independent model reduction via the introduction of the lifted problem structure.
\end{abstract}



\section{Introduction}
\label{sec:intro}

Reduced-order models (ROMs) are an essential enabler for design and optimization of aerospace systems, providing a rapid simulation capability that retains the important dynamics resolved by a more expensive high-fidelity model. Despite a growing number of successes, there remains a tremendous divide between rigorous theory---well developed for the linear case---and the challenging nonlinear problems that are of practical relevance in aerospace applications. For linear systems, ROMs are theoretically well-understood (error analysis, stability, structure preservation) as well as computationally efficient \cite{moore81principal,antoulas05,hesthaven2016certified,volkwein2013model,BGW15surveyMOR}.
For general nonlinear systems, the proper orthogonal decomposition (POD) has been successfully applied to several different problems, but its success typically depends on careful selection of tuning parameters related to the ROM derivation process. For example, nonlinear problems often do not exhibit monotonic improvements in accuracy with increased dimension of the ROM; indeed for some cases, increasing the resolution of the ROM can lead to a numerically unstable model which is practically of no use \cite[Sec.IV.A]{huangAIAA18RomRocketCombustion} as well as \cite{rempfer2000low,bergmann2009enablers} and the references therein. In this paper, we propose an approach to bridge this divide: we show that a general nonlinear system can be transformed into a polynomial form through the process of \emph{lifting}, which introduces auxiliary variables and variable transformations. The lifted system is equivalent to the original nonlinear system, but its polynomial structure offers a number of key advantages.

Ref.~\cite{gu2011qlmor} introduced the idea of lifting nonlinear dynamical systems to quadratic-bilinear (QB) systems for model reduction, and showed that the number of auxiliary variables needed to lift a system to QB form is linear in the number of elementary nonlinear functions in the original state equations.
The idea of variable transformations to promote system structure \add{can be found across different communities, spanning several decades of work}. Ref.~\cite{mccormick1976computability} introduced variable substitutions to solve non-convex optimization problems. \add{Ref.~\cite{kerner1981universal} introduced variable transformations to bring general ordinary differential equations (ODEs) into Riccati form in an attempt to unify theory for differential equations. Ref.~\cite{savageau1987recasting} showed that all ODE systems with (nested) elementary functions can be recast in a special polynomial system form, which are then faster to solve numerically.}
The idea of transforming a general nonlinear system into a system with more structure is also common practice in the control community: the concept of feedback linearization transforms a general nonlinear system into a structured linear model \cite{jakubczyk1980linearization,Khalil_NonlinearSystems}. This is done via a state transformation, where the transformed state might be augmented (i.e., might have increased dimension relative to the original state). However, the lifting transformations known in feedback linearization are specific to the desired model form, and are not applicable in our work here.
In the dynamical systems community, the Koopman operator is a linear infinite dimensional operator that describes the dynamics of observables of nonlinear systems. With the choice of the right observables, linear analysis of the infinite-dimensional Koopman operator helps identify finite dimensional nonlinear state-space dynamics, see Refs.~\cite{rowley09spectralDMD,schmid2010dynamic,mezic13analysisKoopman,tu2013dynamic,kramer17DMD}.

Lifting has been previously considered as a way to obtain QB systems for model reduction in Refs.~\cite{bennerBreiten2015twoSided,bennergoyal2016QBIRKA,bennerGoyal2017BT_quadBilinear}. However, the models considered therein always resulted in a QB system of ordinary differential equations (QB-ODEs), and only one auxiliary lifting variable was needed to yield a QB-ODE. Here, we present a multi-step lifting transformation that leads to a more general class of lifted systems. In particular, for the aerospace example considered in this paper, the system is lifted either to a QB system of differential algebraic equations (QB-DAEs) or to a quartic systems of ODEs.
We then perform POD-based model reduction on this lifted system, exploiting the newly obtained structure.
There are a number of important advantages to reducing a polynomial, and in particular QB, system.
First, ROMs for polynomial systems do not require approximation of the nonlinear function through sampling, since all reduced-order operators can be precomputed. This is in contrast to a general nonlinear system, where an additional approximation step is needed to obtain an efficient ROM~\cite{deim2010,barrault2004empirical,astrid2008missing,grepl2007efficient,carlberg2013gnat,nguyen2008best}.
This property of polynomial ROMs has been exploited in the past, for example, for the incompressible Navier-Stokes equations with quadratic nonlinearities~\cite{holmes_lumley_berkooz_1996,graham1999optimal}, and in the trajectory piecewise linear method~\cite{rewienskiwhite}.
Second, promising progress has been made recently in specialized model reduction for QB systems, such as moment matching~\cite{gu2011qlmor,bennerBreiten2015twoSided}, the iterative rational Krylov algorithm~\cite{bennergoyal2016QBIRKA}, and balanced truncation~\cite{bennerGoyal2017BT_quadBilinear}. The structure of QB systems makes them amenable to input-independent reduced-order modeling, an important feature for control systems and systems that exhibit significant input disturbances.
Third, reducing a structured system is promising in terms of enabling rigorous theoretical analysis of ROM properties, such as stability and error analysis.

\begin{figure}
	\centering
	\includegraphics[width=14cm]{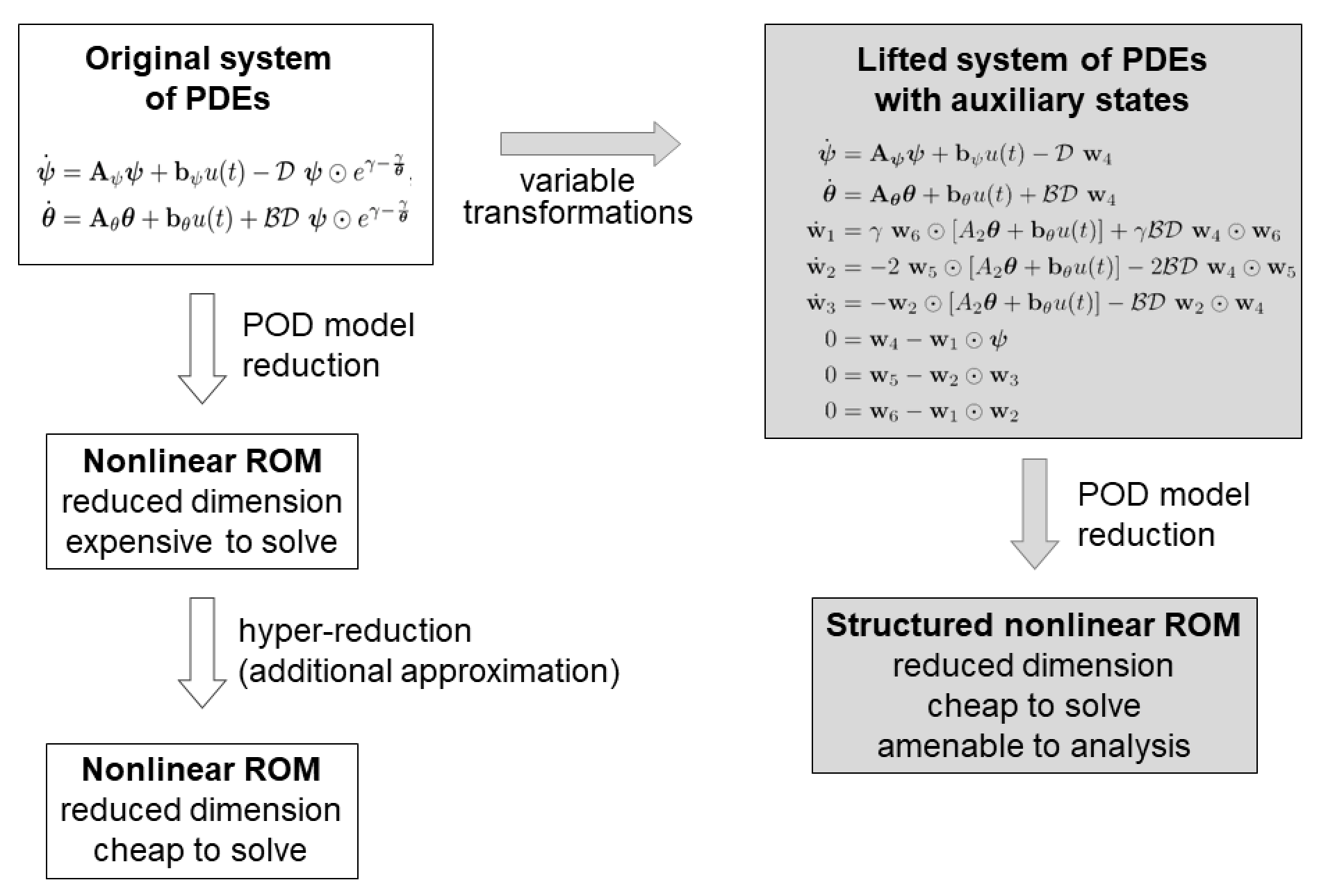}
	\caption{Existing nonlinear model reduction approaches (left flow) require additional approximation of the nonlinear terms; our approach (right flow) first introduces variable transformations to lift the governing equations to a system with more structure, as illustrated for the equations governing the dynamics of a tubular reactor.}
	\label{fig:liftingfigure}
\end{figure}

In this work, our first main contribution is to derive two lifted systems for a strongly nonlinear model of a tubular reactor that models a chemical process. The first lifted model is a quartic ODE. We show that if the goal is to further reduce the polynomial order from quartic to quadratic, then algebraic equations are required to keep the model size of a QB model moderate. Thus, our second lifted model is a QB-DAE. The lifting transformations are nontrivial and proceed in multiple layers.
Our second main contribution is to present a POD-based model reduction method applied to the lifted system. POD is a particularly appropriate choice for the model reduction step (in contrast to previous work which uses balanced truncation and rational Krylov methods), due to the flexibility of the POD approach. In particular, we show that for both the quartic ODE and the QB-DAEs, our POD model reduction method retains the respective structure in the reduction process.
Third, we present numerical comparisons to state-of-the-art methods in nonlinear model reduction. Our lifted ROMs are competitive with state-of-the-art; however, as mentioned above, the structured (polynomial or quadratic) systems have several other advantages.
Figure~\ref{fig:liftingfigure} illustrates our approach and puts it in contrast to state-of-the-art model reduction methods for nonlinear systems.

This paper is structured as follows: Section~\ref{sec:2} briefly reviews POD model reduction, defines polynomial systems and QB-DAEs, and presents the POD-based model reduction of such systems. Section~\ref{sec:lifting} presents the method of lifting general nonlinear systems to polynomial systems, with a particular focus on the case of QB-DAEs.  Section~\ref{sec:FHN} demonstrates and compares the lifting method with state-of-the-art POD-DEIM model reduction for the benchmark problem of the FitzHugh-Nagumo system. Section~\ref{sec:TubReactor} presents the tubular reactor model for which two alternative lifted models are obtained, namely a quartic ODE and a QB-DAE. Numerical results for both cases are compared with POD-DEIM. Finally, Section~\ref{sec:conclusion} concludes the paper.

\section{Polynomial Systems and Proper Orthogonal Decomposition Model Reduction} \label{sec:2}

Section~\ref{sec:POD} briefly reviews the POD method and its challenges. In Section~\ref{sec:polySys} we introduce polynomial systems of ODEs and POD model reduction for such systems. Section~\ref{sec:QB} formally introduces QB-ODE and QB-DAE systems, wich are polynomial systems of order two, but in the latter case with algebraic constraints embedded. That section also presents structure-preserving model reduction for the QB-DAE systems via POD. The quartic, QB-ODE and QB-DAE forms all appear in our applications in  Section~\ref{sec:FHN} and Section~\ref{sec:TubReactor}.

\subsection{Proper Orthogonal Decomposition Model Reduction} \label{sec:POD}

Consider a large-scale nonlinear dynamical system of the form
\begin{equation}
	\dot{\x} = f(\x) + \B \u, \label{eq:DS}
\end{equation}
where $\x(t) \in \mathbb{R}^{n}$ is the state of (large) dimension $n$, $t\geq 0$ denotes time, $\u(t) \in \mathbb{R}^m$ is a time-dependent input of dimension $m$, $\B\in \mathbb{R}^{n\times m}$ is the input matrix,  the nonlinear function $f(\cdot): \mathbb{R}^n \mapsto \mathbb{R}^n$ maps the state $\x$ to $f(\x)$, and $\dot{\x} = \frac{d\x}{dt}$ denotes the time derivative.
Equation~\eqref{eq:DS} is a general form that arises in many engineering contexts. Of particular interest are the systems arising from discretization of partial differential equations. In these cases, the state dimension $n$ is large and simulations of such models are computationally expensive. 
Consequently, we are interested in approximating the full-order model (FOM) in Equation~\eqref{eq:DS} by a ROM of drastically reduced dimension $r\ll n$.

The most common nonlinear model reduction method, proper orthogonal decomposition (POD), computes a basis using snapshot data (i.e., representative state solutions) from simulations of the FOM, see Refs.~\cite{lumley1967structure, sirovich87turbulence, holmes_lumley_berkooz_1996}.
POD has had considerable success in application to aerospace systems (see e.g., \cite{Dowell2001,Lieu2006,Lieu2007, tadmor2007low, BuiThanh2008_AIAA, amsallem2010towards, brunton_rowley_williams_2013,berger2013reduced, taira2017modal}).
Denote the POD basis matrix as $\V\in \mathbb{R}^{n\times r}$, which contains as columns $r$ POD basis vectors. $\V$ is computed from a matrix of $M$ solution snapshots, i.e., $\X = [\x(t_0), \x(t_1), \ldots, \x(t_M)]$. In the case where we have fewer snapshots than states, i.e., $M\ll n$, the simplest form of POD takes the singular value decomposition $\X = \mathbf{U} \mathbf{\Sigma} \mathbf{W}^\top$ and chooses the first $r$ columns of $\mathbf{U}$ to be the POD basis matrix $\V=\mathbf{U}(:,1:r)$. Alternatively, the method of snapshots by Sirovich can be employed~\cite{sirovich87turbulence} to compute $\V$. Regardless, the POD approximation of the state is then
\begin{equation}
	\x\approx \V\widehat{\x},
\end{equation}
where $\widehat{\x}(t) \in \mathbb{R}^{r}$ is the reduced-order state of (small) dimension $r$. Substituting this approximation into Equation~\eqref{eq:DS} and enforcing orthogonality of the resulting residual to the POD basis via a standard Galerkin projection yields the POD ROM
\begin{equation}
	\dot{\widehat{\x}} = \widehat{f}(\widehat{\x}) + \widehat{\B} \u, \label{eq:DSROM}
\end{equation}
with $\widehat{\B} = \V^\top \B \in \mathbb{R}^{r\times m}$, and  $\widehat{f}(\cdot): \mathbb{R}^r \mapsto \mathbb{R}^r$ with $\widehat{f}(\widehat{\x}) =\V^\top f(\V\widehat{\x})$.

Equation~\eqref{eq:DSROM} reveals a well-known challenge with nonlinear model reduction: the evaluation of $\V^\top f(\V\widehat{\x})$ still scales with the FOM dimension $n$. To remedy this problem, state-of-the-art in nonlinear model reduction introduces a second layer of approximation, sometimes referred to as ``hyper-reduction.''
Several nonlinear approximation methods have been proposed, see Refs.~\cite{deim2010,barrault2004empirical,astrid2008missing,grepl2007efficient,carlberg2013gnat,nguyen2008best}, all of which are based on evaluating the nonlinear function $f(\cdot)$ at a sub-selection of sampling points.
Of these, the Discrete Empirical Interpolation Method (DEIM) in Ref.~\cite{deim2010} has been widely used in combination with POD (POD-DEIM), and has been shown to be effective for nonlinear model reduction over a range of applications.
The number of sampling points used in these hyper-reduction methods often scales with the reduced-order model dimension, which leads to an efficient ROM. However, problems with strong nonlinearities can require a high number of sampling points (sometimes approaching the FOM dimension $n$), rendering the nonlinear function evaluations expensive. This has been observed in the case of ROMs for complex flows in rocket combustion engines in Ref.~\cite{huangAIAA18RomRocketCombustion}. A second problem with hyper-reduction is that it introduces an additional layer of approximation to the ROM, which in turn can hinder rigorous analysis of ROM properties such as stability and errors.

\subsection{Polynomial Systems and Proper Orthogonal Decomposition} \label{sec:polySys}

Having discussed nonlinear model reduction via POD in its most general form, we now develop POD models for the specific case of nonlinear systems with polynomial nonlinearities. We will show in Section~\ref{sec:lifting} that lifting transformations can be applied to general nonlinear systems to convert them to this form. We develop here POD models for polynomial systems of order four (quartic systems) and two (quadratic systems), as those arise in our applications; however, the material below extends straightforwardly (at the expense of heavier notation) to the general polynomial case. In the following, the notation $\otimes$ denotes the Kronecker product of matrices or vectors.

A quartic FOM with state $\x(t)$ of dimension $n$ and input $\u(t)$ of dimension $m$  is given by
\begin{align}
	\dot{\x} & = \underbrace{\A \x + \B \u}_{\text{linear}} + \underbrace{\G^{(2)} (\x \otimes \x)}_{\text{quadratic}} + \underbrace{\G^{(3)} (\x\otimes \x \otimes \x)}_{\text{cubic}} + \underbrace{\G^{(4)} (\x\otimes \x \otimes \x \otimes \x)}_{\text{quartic}} \nonumber \\ & + \underbrace{\sum_{k=1}^{m}\N_k^{(1)} \x u_k}_{\text{bilinear}} + \underbrace{\sum_{k=1}^{m}\N_k^{(2)}(\x \otimes \x) u_k}_{\text{quadratic-linear}} , \label{eq:quarticODE}
\end{align}
with $\B \in \mathbb{R}^{n\times m}$ and $\A\in \mathbb{R}^{n\times n}$, and $\G^{(i)}, \N_k^{(i)} \in \mathbb{R}^{n\times n^i}$.  In this form, the matrix $\A$ represents the terms that are linear in the state variables, the matrix $\B$ represents the terms that are linear with respect to the input, the matrices $\G^{(i)}, i=2, \ldots, 4$ represent matricized higher-order tensors for the quadratic, cubic and quartic terms, and the matrices $\N_k^{(1)}$ and $\N_k^{(2)}$ represent respectively the bilinear and quadratic-linear coupling between state and input, with one term for each input $u_k, \ k=1,\ldots, m$.

To reduce the quartic FOM~\eqref{eq:quarticODE}, approximate $\x \approx \V \widehat{\x}$ in the POD basis $\V$ and perform a standard Galerkin projection as described in Section~\ref{sec:POD}, leading to the ROM
\begin{equation}
	\dot{\widehat{\x}} = \widehat{\A} \widehat{\x} + \widehat{\B} \u + \widehat{\G}^{(2)} (\widehat{\x} \otimes \widehat{\x}) + \widehat{\G}^{(3)} (\widehat{\x}\otimes \widehat{\x} \otimes \widehat{\x}) + \widehat{\G}^{(4)} (\widehat{\x}\otimes \widehat{\x} \otimes \widehat{\x} \otimes \widehat{\x}) + \sum_{k=1}^{m}\widehat{\N}_k^{(1)} \widehat{\x} u_k + \sum_{k=1}^{m}\widehat{\N}_k^{(2)}  (\widehat{\x} \otimes \widehat{\x}) u_k. \label{eq:quarticROM}
\end{equation}
The reduced-order matrices and tensors are all straightforward projections of their FOM counterparts onto the POD basis:
$\widehat{\A} = \V^\top \A \V, \ \widehat{\B} = \V^\top \B, \ \widehat{\G}^{(2)} = \V^\top \G^{(2)}(\V \otimes \V), \ \widehat{\G}^{(3)} = \V^\top \G^{(3)}(\V \otimes \V \otimes \V ), \ \widehat{\G}^{(4)} = \V^\top \G^{(4)}(\V \otimes \V \otimes \V \otimes \V)$, $\widehat{\N}_k^{(1)} = \V^\top \N_k^{(1)}\V$, and $\widehat{\N}_k^{(2)} = \V^\top \N_k^{(2)}(\V \otimes \V) $.  Note, that all these reduced-order matrices and tensors can be pre-computed once the POD basis $\V$ is chosen; thus, the POD ROM for the polynomial system recovers an efficient offline-online decomposition and does not require an extra step of hyper-reduction.
\add{Nevertheless, despite Equation~\eqref{eq:quarticROM} preserving the polynomial structure of the original model~\eqref{eq:quarticODE}, the model reduction problem remains challenging. In particular, the training data for POD basis computation, the number of selected modes (especially for problems with multiple variables), and the properties of the model itself (manifested in the system matrices) can all influence the quality of the ROM.}

\subsection{Quadratic-Bilinear Systems and Proper Orthogonal Decomposition} \label{sec:QB}
As a special case of polynomial systems, we focus on quadratic-bilinear (QB) systems for reasons mentioned in Section~\ref{sec:intro}. Consider a system with state $\x(t)$ of dimension $n$ and input $\u(t)$ of dimension $m$.
The general form of a QB system is written
\begin{align}
	\E \dot{\x} & = \underbrace{\A \x + \B \u}_{\text{linear}} + \underbrace{\HH (\x \otimes \x)}_{\text{quadratic}} + \underbrace{\sum_{k=1}^m \N_k \x u_k}_{\text{bilinear}} , \label{eq:qb-dae}
\end{align}
with $ \E \in \mathbb{R}^{n\times n}, \ \A \in \mathbb{R}^{n\times n}, \ \B \in \mathbb{R}^{n\times m}, \ \HH \in \mathbb{R}^{n\times n^2}$ and $\N_k \in \mathbb{R}^{n\times n},\  k=1,\ldots, m$.
The matrices have the same meaning as in the quartic case, except that we use the usual notation $\HH$ for the matricized tensor  that represents the terms that are quadratic in the state variables. In addition, we have introduced the matrix $\E$ (sometimes called the ``mass matrix'') on the left side of the equation.

If the matrix $\E$ is nonsingular, then Equation~\eqref{eq:qb-dae} is a QB system of ODEs. If the matrix $\E$ is singular, then Equation~\eqref{eq:qb-dae} is a QB system of differential algebraic equations (DAEs)\footnote{Note that when the system is a DAE, $\x(t)$ is not technically a ``state'' in the sense of being the smallest possible number of variables needed to represent the system; however, it is common in the literature to still refer to $\x(t)$ as the ``state'', as we will do here.}; in particular $\E$ will have zero rows corresponding to any algebraic equations.

We now focus on the QB-DAE case, as such a system arises from lifting transformations, as we see later for the tubular reactor model in Section~\ref{sec:TubQBDAE}. The QB-DAE state is partitioned as $\x = [\x_1^\top \ \x_2^\top]^\top$ with $\x_1 \in \mathbb{R}^{n_1}$ being the dynamically evolving states and $\x_2 \in \mathbb{R}^{n_2}$ the algebraically constrained variables, with $n=n_1+n_2$. A lifting transformation resulting in QB-DAEs often leads to matrices with special structure as follows:
\begin{equation}
	\E  = \begin{bmatrix} \E_{11} & \bzero \\ \bzero& \bzero\end{bmatrix}, \quad
	\A  = \begin{bmatrix} \A_{11} & \A_{12}\\  \bzero & \I_{n_2} \end{bmatrix}, \quad
	\HH = \begin{bmatrix} \HH_1 \\ \HH_{2} \end{bmatrix} , \quad
	\N_k = \begin{bmatrix} \N_{k,11} & \N_{k,12} \\ \bzero & \bzero \end{bmatrix}, \quad
	\B  = \begin{bmatrix} \B_1 \\ \bzero \end{bmatrix}. \label{eq:QBDAE-matrices}
\end{equation}
Here, $\I_{n_2}$ is the $n_2\times n_2$ identity matrix and $\bzero$ denotes a matrix of zeros of appropriate dimension. Moreover, $\B_1 \in \mathbb{R}^{n_1 \times m} $ and $\A_{11}, \E_{11}, \N_{11} \in \mathbb{R}^{n_1\times n_1}$.
The QB-DAE with the above structure can then be rewritten as
\begin{align}
	\E_{11} \dot{\x }_1 &= \A_{11} \x _1 + \A_{12} \x _2 + \B_1 \u + \HH_1 (\x \otimes \x ) + \sum_{k=1}^m \N_{k,11} \x _1u_k + \N_{k,12} \x _2u_k  \label{eq:ODE_x1}, \\
	\bzero & = \x _2 - \widetilde{\HH}_2 (\x_1\otimes\x_1)\label{eq:DAE_x2} ,
\end{align}
where $\widetilde{\HH}_2 \in \mathbb{R}^{n_2\times n_2^2}$ is obtained from $\HH_2 \in \mathbb{R}^{n_2\times n^2} $ by deleting columns corresponding to the zeros in the Kronecker product. We note that Equation~\eqref{eq:ODE_x1} is the $n_1$th-order system of ODEs describing dynamical evolution of the states $\x_1$, while Equation~\eqref{eq:DAE_x2} are the $n_2$ algebraic equations that enforce the relationship between the constrained variables $\x_2$ and the states $\x_1$.
%

The QB-DAE~\eqref{eq:qb-dae}--\eqref{eq:QBDAE-matrices} can be directly reduced using a POD projection.
To retain the DAE structure in the model, we use the projection matrix
\begin{equation} \label{eq:defV}
	\V =  \begin{bmatrix} \V_1 & \bzero  \\ \bzero & \V_2  \end{bmatrix},
\end{equation}
where $\V_1 \in \mathbb{R}^{n_1\times r_2}$ and $\V_2 \in \mathbb{R}^{n_2\times r_2}$ are the POD basis matrices that contain as columns POD basis vectors for $\x_1$ and $\x_2$, respectively, and $r_1+r_2=r$.
We approximate the state $\x \approx \V\widehat{\x}$ where $\widehat{\x} \in \mathbb{R}^r$ is the reduced state of dimension $r \ll n$. By definition, $\x_1 \approx \V_1\widehat{\x}_1$ and $\x_2 \approx \V_2\widehat{\x}_2$.
Introducing this approximation to (\ref{eq:qb-dae}) and using the standard POD Galerkin projection yields the reduced-order model
\begin{align} \label{eq:QBROM}
	\widehat{\E} \dot{\widehat{\x}} & = \widehat{\A} \widehat{\x} + \widehat{\B} \u + \widehat{\HH} (\widehat{\x} \otimes \widehat{\x}) + \sum_{k=1}^m \widehat{\N}_k \widehat{\x} u_k .
\end{align}
The reduced-order matrices can be pre-computed as
\begin{equation*}
	\widehat{\E}  = \begin{bmatrix} \widehat{\E}_{11} & \bzero \\ \bzero& \bzero\end{bmatrix}, \quad
	\widehat{\A}  = \begin{bmatrix} \widehat{\A}_{11} & \widehat{\A}_{12}\\  \bzero & \I_{r_2} \end{bmatrix}, \quad
	\widehat{\HH} = \begin{bmatrix} \widehat{\HH}_1 \\ \widehat{\HH}_{2} \end{bmatrix} , \quad
	\widehat{\N}_k = \begin{bmatrix}  \widehat{\N}_{k,11} & \widehat{\N}_{k,12} \\ \bzero & \bzero \end{bmatrix}, \quad
	\widehat{\B}  = \begin{bmatrix} \widehat{\B}_1 \\ \bzero  \end{bmatrix},
\end{equation*}
where $ \widehat{\E}_{11} = \V_1^\top \E_{11} \V_1, \ \widehat{\A}_{11} = \V_1^\top\A_{11} \V_1, \ \widehat{\A}_{12}= \V_1^\top\A_{12} \V_2, \  \widehat{\N}_{k,11} = \V_1^\top \N_{k,11} \V_1,$ \\ $\widehat{\N}_{k,12} = \V_1^\top \N_{k,12} \V_2 , \  \widehat{\B}_1 = \V_1^\top \B_1$.
The quadratic tensors can be precomputed as
\begin{align}
	\widehat{\HH}_1 & = \V_1^\top \HH_1 \left (\begin{bmatrix} \V_1 & \bzero  \\ \bzero & \V_2  \end{bmatrix} \otimes \begin{bmatrix} \V_1 & \bzero  \\ \bzero & \V_2  \end{bmatrix} \right ) \in \mathbb{R}^{r_1 \times (r_1 + r_2)^2 }, \\
	\widehat{\HH}_2 & = \V_2^\top \widetilde{\HH}_{2} (\V_1 \otimes \V_1) \in \mathbb{R}^{r_2\times r_1^2 }.
\end{align}
The ROM can then be rewritten as
\begin{align}
	\widehat{\E}_{11} \dot{\widehat{\x} }_1 & = \widehat{\A}_{11} \widehat{\x} _1 + \widehat{\A}_{12} \widehat{\x}_2 + \widehat{\B}_1 \u + \widehat{\HH}_1 \left (\begin{bmatrix}
		\widehat{\x}_1 \\ \widehat{\x}_2 \end{bmatrix} \otimes \begin{bmatrix}
		\widehat{\x}_1 \\ \widehat{\x}_2 \end{bmatrix} \right ) + \sum_{k=1}^m \widehat{\N}_{k,11} \widehat{\x}_1 u_k + \widehat{\N}_{k,12} \widehat{\x}_2u_k  \label{eq:ROM_x1},\\
	\bzero & = \widehat{\x} _2 - \widehat{\HH}_2 (\widehat{\x}_1\otimes \widehat{\x}_1). \label{eq:ROM_x2}
\end{align}
With this projection, the index of the DAE is preserved, since the structure of the algebraic equations remains unaltered.
Since all ROM matrices and tensors can be precomputed, no additional approximations (e.g., DEIM, other hyper-reduction) are needed.
The solution of this system is described in the Appendix. Note that as a special case, if $\V_2=\I$ we can obtain a quartic ROM by eliminating the algebraic constraint and inserting $\widehat{\x}_2$ ($=\x_2$) from Equation~\eqref{eq:ROM_x2} into Equation~\eqref{eq:ROM_x1}.

Having formally introduced QB systems, the next section shows the lifting method applied to nonlinear systems, and how QB systems (DAEs and ODEs) can be obtained in the process.

\section{Lifting Transformations} \label{sec:lifting}

With the formal definition of polynomial and QB systems at hand, we now introduce the concept of lifting and give an example that illustrates the approach.
Lifting is a process that transforms a nonlinear dynamical system with $n$ variables into an equivalent system of $\tilde{n}>n$ variables by introducing $\tilde{n}-n$ additional auxiliary variables. The lifted system has larger dimension, but has more structure. For more details on lifting, we refer the reader to Ref.~\cite{gu2011qlmor}.
Our goal is to transform the original nonlinear model into an equivalent polynomial system via lifting. We target this specific structure, since a large class of nonlinear systems can be written in this form, and since polynomial systems---and as a special case QB systems---are directly amenable to model reduction via POD. Moreover, as illustrated below, lifting to a system of DAEs, instead of requiring the lifted model to be an ODE, keeps the number of auxiliary variables to a manageable level.

The method is best understood with an example.
\begin{example}
	Consider the ODE
	\begin{equation}
		\dot{x}= x^4 + u, \label{eq:ExampleODE}
	\end{equation}
	where $u(t)$ is an input function and $x(t)$ is the one-dimensional state variable. We choose the auxiliary state $w_1=x^2$, which makes the original dynamics~\eqref{eq:ExampleODE} quadratic. The auxiliary state dynamics are (according to the chain rule, or Lie derivative) $\dot{w}_1 = 2x\dot{x} = 2x[w_1^2 + u]$, and hence cubic in the new state $[x, w_1]$. Now, introduce another auxiliary state $w_2 = w_1^2$. Then we have $\dot{w}_1 =  2x[w_1^2 + u] = 2x[w_2+u]$ and $\dot{x}= w_2 + u$. However, we have that $\dot{w}_2 = 2w_1\dot{w}_1 = 4xw_1[w_2 + u]$, which is still cubic. Choosing one additional auxiliary state $w_3 = x w_1$ then makes the overall system QB, since we have $\dot{w}_3 = \dot{x} w_1 + x \dot{w}_1 = [w_2 + u]w_1 + x[2xw_2 + 2xu] = w_1 w_2 + w_1 u + 2w_1w_2 + 2w_1u$. Overall, the nonlinear equation~\eqref{eq:ExampleODE} with one state variable is equivalent to the QB-ODE with four state variables
	\begin{align}
		\dot{x}     & =   w_2  +  u, \label{eq:ex1-1}\\
		\dot{w}_1 & =  2xw_2  +  2xu, \\
		\dot{w}_2 & =  4w_2 w_3  +  4 w_3 u, \\
		\dot{w}_3 & =  3w_1 w_2  + 3w_1 u. \label{eq:ex1-4}
	\end{align}

	An alternative approach is to include the algebraic constraint $ w_1 = x^2$ and instead obtain a QB differential algebraic equation (QB-DAE) with two variables as
	\begin{align}
		\dot{x}     & =   w_1^2  +  u, \label{eq:ex1-dae1}\\
		0 & =  w_1 - x^2. \label{eq:ex1-dae2}
	\end{align}
	We emphasize that the system (\ref{eq:ex1-1})--(\ref{eq:ex1-4}) and the system (\ref{eq:ex1-dae1})--(\ref{eq:ex1-dae2}) are both equivalent to the original nonlinear equation (\ref{eq:ExampleODE}), in the sense that all three systems yield the same solution $x(t)$.
	
	This example illustrates an interesting point in lifting dynamic equations. Even when lifting to a QB-ODE might be possible, our approach of permitting DAEs keeps the number of auxiliary variables low. In particular, Gu~\cite{gu2011qlmor} showed favorable upper bounds for auxiliary variables for lifting to QB-DAEs versus QB-ODEs. This will become important when we consider systems arising from discretization of PDEs, where the number of state variables is already large.
	
\end{example}
The lifted representation is not unique, and we are not aware of an algorithm that finds the \textit{minimal} polynomial system that is equivalent to the original nonlinear system. Moreover, different lifting choices can influence system properties, such as stiffness of the differential equations.

\begin{example}
	Writing the system (\ref{eq:ex1-1})--(\ref{eq:ex1-4}) in the form (\ref{eq:qb-dae}) with $\x=[x \ w_1 \ w_2 \ w_3]^\top $ and the quadratic term $\x\otimes \x = [ x^2 \ xw_1 \ xw_2 \ xw_3 \ w_1 x \ w_1^2 \ w_1 w_2 \ w_1 w_3\ w_2x \ w_2 w_1 \ w_2^2 \ w_2 w_3 $ \\ $\ldots  w_3 x \ w_3 w_1 \ w_3 w_2 \ w_3^2 ]^\top$ yields
	$$
	\E = \left[
	\begin{array}{cccc}
	1 & 0 & 0 & 0 \\
	0 & 1 & 0 & 0 \\
	0 & 0 & 1 & 0 \\
	0 & 0 & 0 & 1 \\
	\end{array}
	\right], \
	\A = \left[
	\begin{array}{cccc}
	0 & 0 & 1 & 0 \\
	0 & 0 & 0 & 0 \\
	0 & 0 & 0 & 0 \\
	0 & 0 & 0 & 0 \\
	\end{array}
	\right], \
	\N_1 = \left[
	\begin{array}{cccc}
	0 & 0 & 0 & 0 \\
	2 & 0 & 0 & 0 \\
	0 & 0 & 0 & 4 \\
	0 & 3 & 0 & 0 \\
	\end{array}
	\right], \
	\B = \left[
	\begin{array}{c}
	1 \\
	0 \\
	0 \\
	0 \\
	\end{array}
	\right],
	$$
	and for the quadratic tensor $\HH\in \mathbb{R}^{4\times 16}$ we have
	$$
	\HH_{2,3} =2, \quad \HH_{3,12} =4, \quad \HH_{4,7} =3, \quad \HH_{i,j} = 0 \ \  \text{otherwise}.
	$$
	Note that this is a system of ODEs (the matrix $\E$ is full rank). In contrast, the system (\ref{eq:ex1-dae1})--(\ref{eq:ex1-dae2}) with $\x=[x \ w_1]^\top$ and $\x \otimes \x = [x^2 \ x w_1 \ w_1 x \ w_1^2]^\top$ yields the DAEs, also of the form (\ref{eq:qb-dae}) but with smaller dimension and singular $\E$, as follows:
	$$
	\E = \left[
	\begin{array}{cc}
	1 &  0 \\
	0 & 0 \\
	\end{array}
	\right], \
	\A = \left[
	\begin{array}{cc}
	0 & 0 \\
	0 & 1 \\
	\end{array}
	\right], \
	\HH = \left[
	\begin{array}{cccc}
	0 & 0 & 0 & 1\\
	-1 & 0 & 0 & 0 \\
	\end{array}
	\right], \
	\N_1 = \left[
	\begin{array}{cc}
	0 & 0 \\
	0 & 0 \\
	\end{array}
	\right], \
	\B = \left[
	\begin{array}{c}
	1 \\
	0 \\
	\end{array}
	\right].
	$$
	Again note that both of these representations are equivalent to the original system~\eqref{eq:ExampleODE}, with no approximation introduced.
	
\end{example}

\section{Benchmark Problem: FitzHugh-Nagumo}
\label{sec:FHN}

This section illustrates our nonlinear model reduction approach on the FitzHugh-Nagumo system, which is a model for the activation and deactivation of a spiking neuron. It is a benchmark model in nonlinear reduced-order modeling, and has been explored in the context of DEIM in Ref.~\cite{deim2010}, balanced model reduction in Ref.~\cite{bennerGoyal2017BT_quadBilinear}, and interpolation-based model reduction in Ref.~\cite{bennerBreiten2015twoSided}.

\subsection{FitzHugh-Nagumo Problem Definition}
The FitzHugh-Nagumo governing PDEs are
\begin{align}
	\epsilon \dot{v} & = \epsilon^2 v_{ss}  - v^3 + 0.1 v^2 - 0.1v - w + c, \label{eq:FHN1}\\
	\dot{w} &  = hv - \gamma w + c, \label{eq:FHN2}
\end{align}
where $s\in [0,L]$ is the spatial variable and the time horizon of interest is $t \in [0,t_f]$. The states of the system are voltage $v(s,t)$ and recovery of voltage $w(s,t)$. The notation $v_{ss}(s,t) := \frac{\partial^2}{\partial s^2}v(s,t)$ denotes a second order spatial derivative; similarly, $v_s(s,t)$ denotes a first spatial derivative.
The initial conditions are specified as
\begin{align*}
	v(s,0)    & =0, 									 	 & w(s,0)  =0,  & \quad s\in [0,L],
\end{align*}
and the boundary conditions are
\begin{align*}
	v_s(0,t) & = u(t),  & v_s(L,t)  =0, & \quad t\geq 0,
\end{align*}
where $u(t) = 5\times 10^4 \ t^3 \exp (-15t)$.  In the problem setup we consider, the parameters are given by \add{$L=1$}, $c=0.05$, $\gamma =2$, $h=0.5$, and $\epsilon = 0.015$.
%

\subsection{FitzHugh-Nagumo Lifted Formulation}
To lift the FitzHugh-Nagumo equations to QB form, we follow the same intuitive lifting as in Ref.~\cite{bennerBreiten2015twoSided}. Choose $z=v^2$, which renders the original equations~\eqref{eq:FHN1}--\eqref{eq:FHN2} quadratic.
The auxiliary equation becomes
\begin{align*}
	\dot{z} & = 2 v \dot{v}  = 2 v[\epsilon^2 v_{ss}  - v^3 + 0.1 v^2 - 0.1v - w + c] \\
	& = 2 [\epsilon^2 vv_{ss}  - z^2 + 0.1 z v - 0.1 z - wv+ cv],
\end{align*}
and is quadratic in the new variable.
The lifted QB system then reads as
\begin{align*}
	\epsilon \dot{v} & = \epsilon^2 v_{ss}  - z v + 0.1 z - 0.1v - w + c,\\
	\dot{w} &  = hv - \gamma w + c, \\
	\dot{z} & =2 [\epsilon^2 vv_{ss}  - z^2 + 0.1 z v - 0.1 z  - w v + cv].
\end{align*}
The initial conditions for the auxiliary variable need to be consistent, i.e., $z(s,0) = v(s,0)^2, s \in [0,L]$.
The boundary conditions are obtained by applying the chain rule:
$$	
z_s(L,t) = 2 v(L,t) \underbrace{v_s(L,t)}_{=0} = 0,
$$	
and on the left side
\begin{equation*}
	z_s(0,t) = 2 v(0,t) \underbrace{v_s(0,t)}_{u(t)} = 2 v(0,t)u(t). \label{eq:BC}
\end{equation*}

The full model is discretized using finite differences, where each variable is discretized with \add{$n=512$} degrees of freedom, i.e., the overall dimension of the QB model is $3n=1536$. The resulting QB-ODE system is
\begin{equation*}
	\E \dot{\x} = \A \x + \B \u + \HH (\x \otimes \x) + \sum_{k=1}^2 \N_k \x u_k ,
\end{equation*}
where $\E=\epsilon \I_{3n}$ is diagonal, $\A,\N_1, \N_2 \in \mathbb{R}^{3n\times3n}$ and $\HH \in \mathbb{R}^{3n\times (3n)^2}$. The input matrix is $\B\in \mathbb{R}^{3n\times 2}$, with the second column of $\B$ being copies of $c$ (the constant in Equations~\eqref{eq:FHN1} and \eqref{eq:FHN2}) and the first column of $\B$ having a 1 at the first entry. Thus, the input $\u = [u(t), 1]$. This benchmark model is available at Ref.~\cite{morwiki_modFHN}.
%

\subsection{FitzHugh-Nagumo Lifted Quadratic-Bilinear Reduced-Order Model}
\begin{figure}[H]
	\begin{subfigure}{}
		\includegraphics[width=7cm]{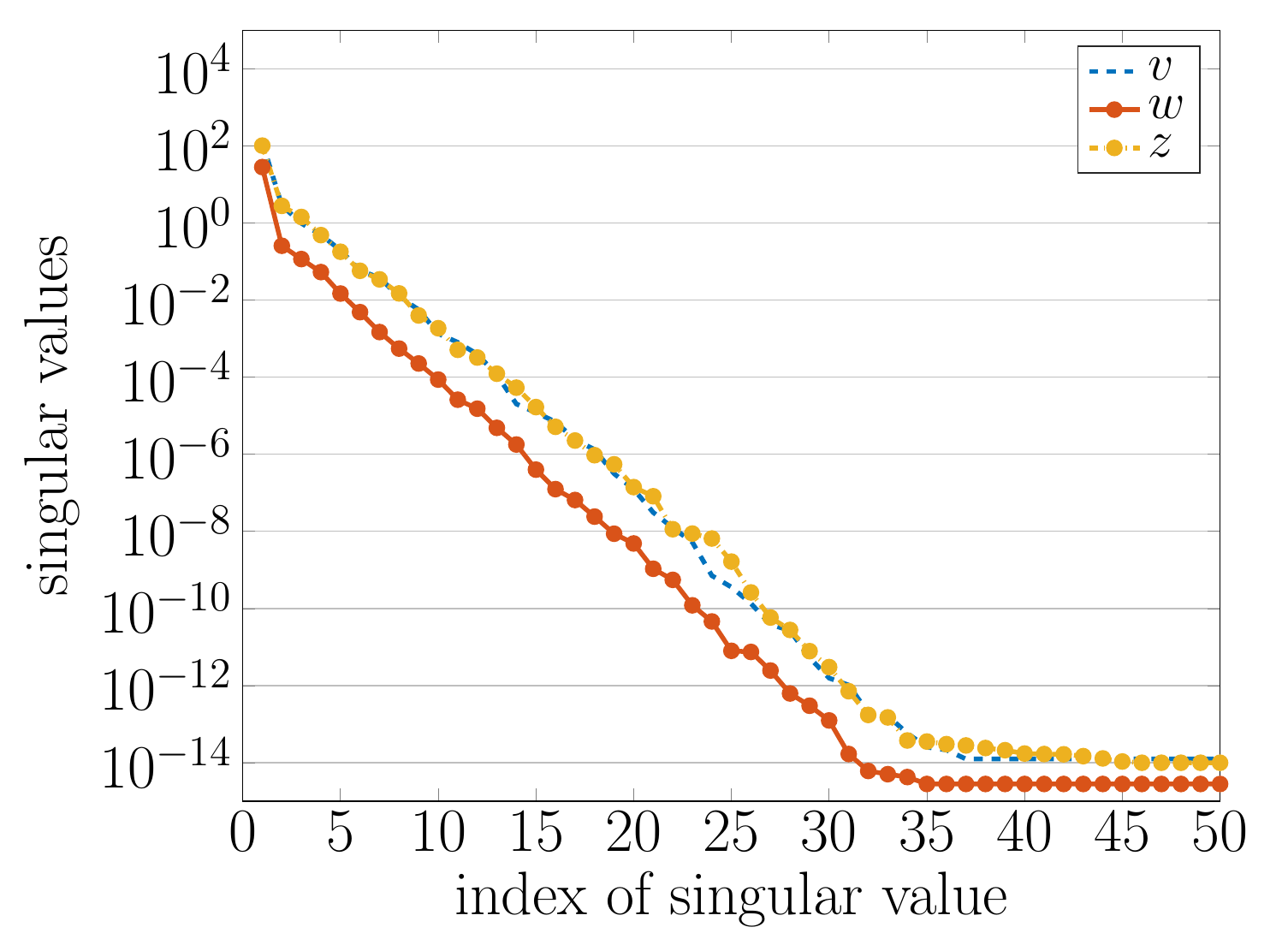}
	\end{subfigure}	
	\begin{subfigure}{}
		\includegraphics[width=7cm]{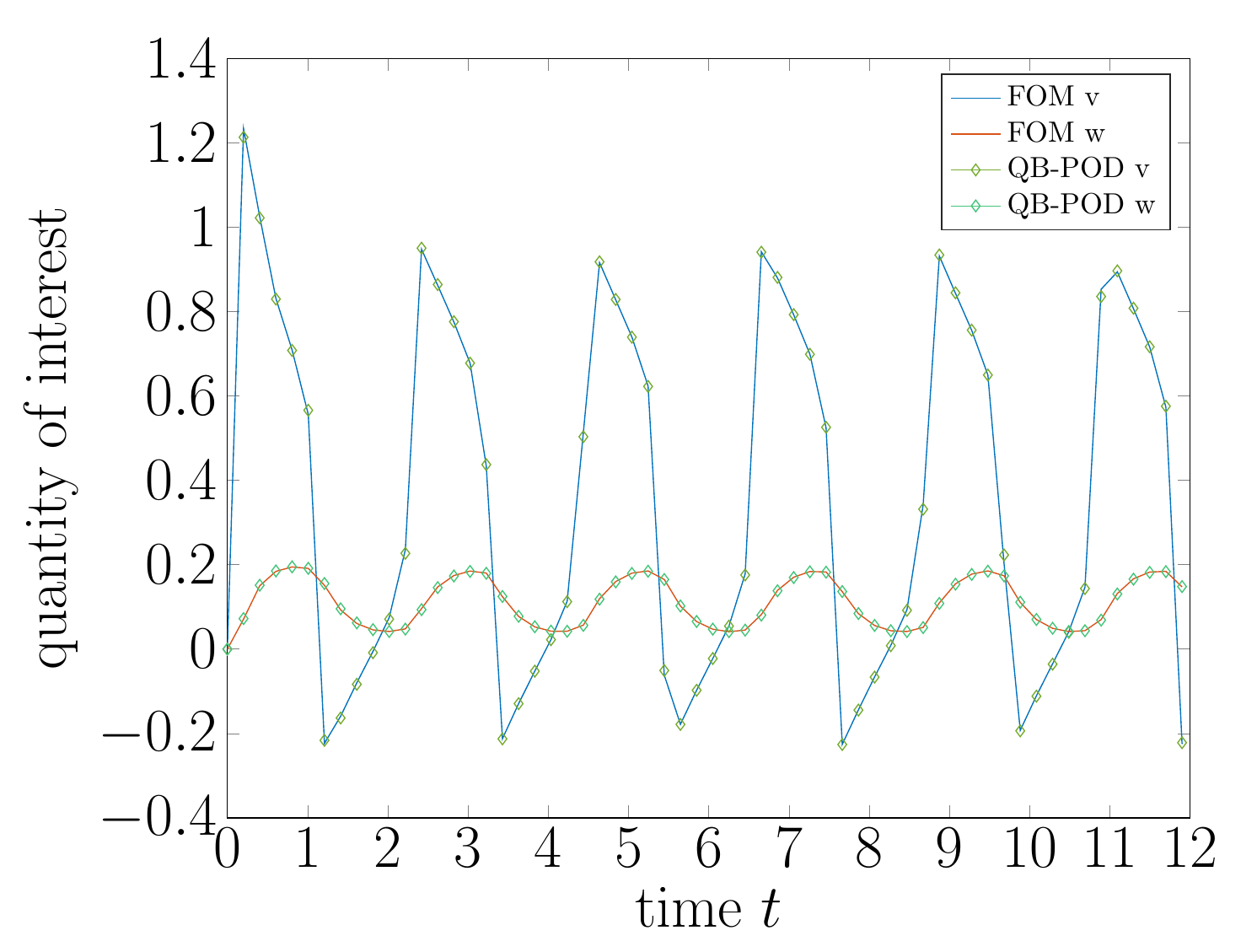}
	\end{subfigure}			
	\caption{FitzHugh-Nagumo system. Left: Decay of singular values of snapshot matrices for three variables. Right: Quantities of interest $w(0,t)$ and $v(0,t)$ comparing FOM simulations and the QB-POD reduced model of dimension $3r=9$.  }
	\label{fig:FHN_SVD_SIM}
\end{figure}

We simulate this lifted full-order system for \add{$t_f = 12s$ and collect $n_t=150$} snapshots of the state solutions at equally spaced times. \add{For the computation of the POD basis, we only use the first 100 snapshots until $t=8s$. Thus, all the ROMs in this section predict 50\% further past the training data}. We compute a separate POD basis for each state variable. This means that for the original system, we compute a POD basis for $v$ and a POD basis for $w$; for the lifted system, we also compute a POD basis for the auxiliary variable $z=v^2$.
Figure~\ref{fig:FHN_SVD_SIM}, left, shows the decay of the singular values for the snapshot matrices of the three state variables, $v$, $w$ and $z$. As expected, the singular values for the snapshot matrix of the auxiliary variable $z=v^2$ show the same decay (up to numerical accuracy) as those for the original variable $v$.
We compute the POD reduced model of the lifted QB system, as described in Section~\ref{sec:QB}. Figure~\ref{fig:FHN_SVD_SIM}, right, shows quantities of interest, namely $w(0,t)$ and $v(0,t)$, computed using the FOM and the QB-POD reduced model with $3r=9$ ($r$ reduced states per variable).
The reduced model captures the limit-cycle oscillations well and is visually indistinguishable from the FOM solution.

\begin{figure}
	\begin{center}
		\includegraphics[width=10cm]{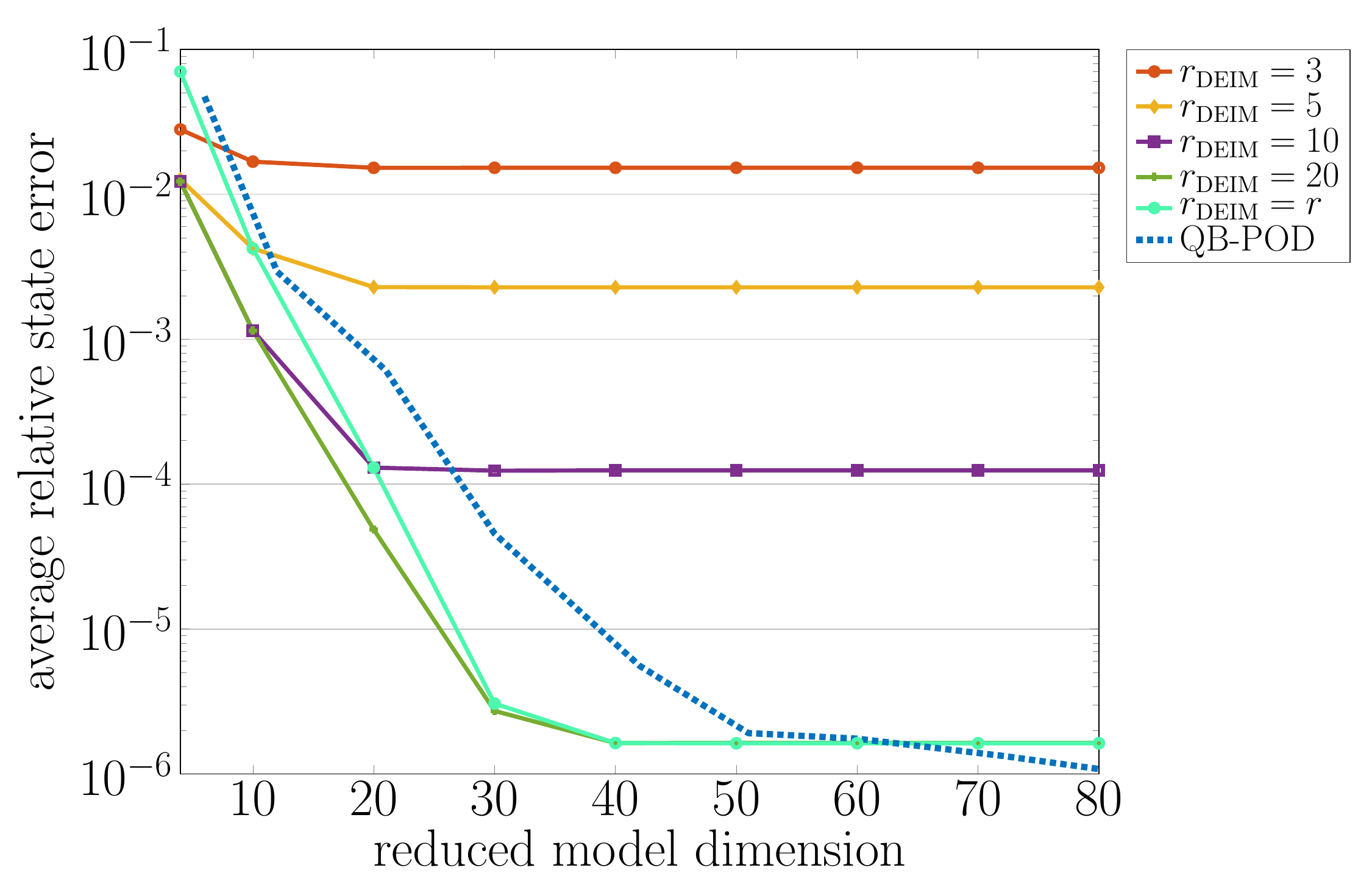}
	\end{center}
	\caption{Average relative state error $\frac{1}{n_t} \sum_{i=1}^{n_t} \Vert \x(t_i)  - \x^{\text{ROM}}(t_i) \Vert / \Vert \x(t_i) \Vert$ as a function of the ROM dimension, for POD-DEIM models and the lifted QB-POD model.}
	\label{fig:FHN_err}
\end{figure}

Figure~\ref{fig:FHN_err} compares the accuracy of the lifted QB-POD models with POD-DEIM models. As first introduced in Ref.~\cite{deim2010}, the POD-DEIM approach reduces the original system with an additional approximation (via DEIM) of the nonlinear term. This approximation is necessary in order for the reduced model to be computationally efficient\footnote{Although in fact we note the cubic nature of the original model, which could directly be exploited---this seems to be overlooked in the literature.}. This requires the following additional steps: (1)~during the full model simulation collect snapshots of the nonlinear term, in addition to snapshots of the states; (2)~apply the POD to the nonlinear term snapshot set to compute the DEIM basis; (3)~select $r_\text{DEIM}$ DEIM interpolation points; and (4)~approximate the projected nonlinear term using the corresponding first $r_\text{DEIM}$ basis vectors. As in Ref.~\cite{deim2010}, we approximate each variable with $r$ basis functions, so the POD-DEIM model has $2r$ dimensions.
Let $\x(t) = [\mathbf{v}(t)^\top, \mathbf{w}(t)^\top]^\top$ be the state vector of the FOM, and $\x^{\text{ROM}}(t)$ be the approximation of that state computed by the different ROM simulations (i.e., $\x^{\text{ROM}}(t)$ contains those components of $\V\widehat{\x}(t)$ that correspond to the original states; for QB-POD we do not measure the error in approximations of the auxiliary variables in order to provide an appropriate comparison). Plotted in Figure~\ref{fig:FHN_err} are the relative errors in the state vector averaged over time, i.e., $\frac{1}{n_t} \sum_{i=1}^{n_t} \Vert \x(t_i)  - \x^{\text{ROM}}(t_i) \Vert / \Vert \x(t_i) \Vert$.  The x-axis plots the overall dimension of the ROM, i.e., the total number of basis functions for the two (for POD-DEIM model) or three (for QB-POD model) variables. For the POD-DEIM models, we show several choices of $r_{\text{DEIM}}$, the number of DEIM interpolation points.
Figure~\ref{fig:FHN_err} shows the characteristic POD-DEIM reduced model error behavior where the number of DEIM interpolation points limits the accuracy of the reduced model, and thus the errors flatten out once a threshold number of POD basis functions is reached. The quality of the reduced model can then only be improved by increasing the number of DEIM interpolation points, which reduces the error in the approximation of the nonlinear term. \add{We also show a POD-DEIM model that increases the DEIM interpolation points with the reduced-dimension, i.e., we have $r_{\text{DEIM}}=r$, yet this model also levels out around $r=35$, so increasing the DEIM dimension does not further improve the model. This is a feature of the FitzHugh-Nagumo problem, as the singular values of the states (see Figure~\ref{fig:FHN_SVD_SIM}) and also the nonlinear snapshots decay to machine precision around $r=35$, see also~\cite{deim2010}}. In contrast, our lifted QB-POD reduced model has no additional approximation step and its error steadily decreases as the number of POD basis functions is increased. These results show that our lifted POD approach recovers the accuracy of a regular POD approach, but has the added benefit that it does not require additional approximation to handle the nonlinear terms.
%

\section{Application: Tubular Reactor Model} \label{sec:TubReactor}

Section~\ref{sec:Tubmodel} describes a tubular reactor model that has strong nonlinearities and limit cycle oscillations representative of those in combustion engines. We demonstrate the benefits of lifting and POD on this problem. First, we bring the system into polynomial form, namely a fourth-order ODE, see Section~\ref{sec:TubquarticODE}. We further lift the polynomial system to a QB-DAE in Section~\ref{sec:TubQBDAE}. Section~\ref{sec:TubPODdetails} presents details for the computation of POD reduced models, and Section~\ref{sec:numerics} shows and discusses the numerical results.

\subsection{Partial Differential Equation Model and Discretization} \label{sec:Tubmodel}
A one-dimensional non-adiabatic tubular reactor model with single reaction is modeled following Refs.~\cite{zhou2012thesis,heinemann1981tubular} as	
\begin{align}
	\dot{\spe}   & = \frac{1}{Pe} \spe_{ss} - \spe_s - \Dam  f(\spe, \theta; \gamma), \label{eq:PDEa} \\
	\dot{\theta} & = \frac{1}{Pe} \theta_{ss} - \theta_s - \beta (\theta - \theta_{\text{ref}}) + \mathcal{B}\Dam f(\spe,\theta; \gamma), \label{eq:PDEb}
\end{align}
with spatial variable $s\in (0,1)$, time $t >0$ and the nonlinear term
\begin{equation*}
	f(\spe, \theta; \gamma) =\spe e^{\gamma - \frac{\gamma}{\theta}}.
\end{equation*}
Robin BCs are imposed on left boundary
\begin{equation*}
	\spe_s(0,t) = Pe (\spe(0,t)-\mu), \qquad \theta_s(0,t) = Pe (\theta(0,t)-1)
\end{equation*}
and Neumann boundary conditions on the right
\begin{equation*}
	\spe_s(1,t) = 0, \quad \theta_s(1,t) =0.
\end{equation*}
The initial conditions are prescribed as
\begin{equation*}
	\spe(s,0) = \spe_0(s), \qquad \theta(s,0) = \theta_0(s).
\end{equation*}	
The variables of the model are the species concentration $\spe$ and temperature $\theta$. The parameters are the Damk{\"o}hler number $\Dam$, P{\`e}clet number $Pe$ as well as known constants $\mathcal{B},\ \beta, \ \theta_{\text{ref}}, \gamma$.
%
%
It is shown in Ref.~\cite{heinemann1981tubular} that when $Pe = 5, \gamma = 25, \mathcal{B} = 0.5, \beta = 2.5, \theta_{\text{ref}} \equiv 1$, the system exhibits a Hopf bifurcation with respect to $\Dam$ in the range $\Dam \in [0.16, 0.17]$; that is, there exists a critical Damk{\"o}hler number $\Dam^{c} = 0.165$ such that for $\Dam^{c} < \Dam$ the unsteady solution eventually converges to a non-trivial steady state.

We discretize the model via finite differences, for details see Ref.~\cite{zhou2012thesis}. The discretized variables are $\boldsymbol{\spe}\in \mathbb{R}^n$ and  $\boldsymbol{\theta}\in \mathbb{R}^n$, so that the resulting dimension of the discretized system is $2n$. The resulting FOM reads as:
\begin{align}
	\dot{\boldsymbol{\spe}} & = \A_\spe \boldsymbol{\spe} + \mathbf{b}_\spe u(t)- \mathcal{D} \ \boldsymbol{\spe} \odot e^{\gamma - \frac{\gamma}{\boldsymbol{\theta}}} , \label{eq:TubODE1} \\
	\dot{\boldsymbol{\theta}} & = \A_\theta \boldsymbol{\theta} + \mathbf{b}_\theta u(t) + \mathcal{B}\mathcal{D} \ \boldsymbol{\spe} \odot e^{\gamma - \frac{\gamma}{\boldsymbol{\theta}}}, \label{eq:TubODE2}
\end{align}
where $\A_{\boldsymbol{\spe}}$ and $\A_{\boldsymbol{\theta}}$ are $n\times n$ matrices, and $\mathbf{b}_\spe, \mathbf{b}_\theta \in \mathbb{R}^n$ enforce the boundary conditions via $u(t)\equiv 1$. Here, we use the (Hadamard) componentwise product of two vectors, i.e., $[\boldsymbol{\spe} \odot\boldsymbol{\theta}]_i = \boldsymbol{\spe}_i\boldsymbol{\theta}_i$. Note that with the exponential nonlinearity, this is a general nonlinear FOM of the form Equation~\eqref{eq:DS}. Direct POD of this model would require additional approximation of the nonlinear term (e.g., via DEIM).

\subsection{Lifted Model 1: A Quartic Ordinary Differential Equation} \label{sec:TubquarticODE}
We start with polynomializing the system via the dependent variables
\begin{align} \label{eq:wDef}
	\mathbf{w}_1 = e^{\gamma - \frac{\gamma}{\boldsymbol{\theta}}}, \quad \mathbf{w}_2 = \boldsymbol{\theta}^{-2}, \quad \mathbf{w}_3 = \boldsymbol{\theta}^{-1}.
\end{align}
Application of the chain rule yields
\begin{align*}
	\dot{\mathbf{w}}_1 & = \mathbf{w}_1\odot (\gamma \ \boldsymbol{\theta}^{-2})\odot \dot{\boldsymbol{\theta}} = \gamma \mathbf{w}_1\odot \mathbf{w}_2 \odot\dot{\boldsymbol{\theta}}, \\
	\dot{\mathbf{w}}_2 & = -2 \boldsymbol{\theta}^{-3}\odot \dot{\boldsymbol{\theta}} = -2 \mathbf{w}_2\odot \mathbf{w}_3\odot \dot{\boldsymbol{\theta}}, \\
	\dot{\mathbf{w}}_3 & = - \boldsymbol{\theta}^{-2}\odot \dot{\boldsymbol{\theta}} = - \mathbf{w}_2\odot \dot{\boldsymbol{\theta}},
\end{align*}
where $\dot{\boldsymbol{\theta}}$ is given by the right-hand-side of Equation~\eqref{eq:TubODE2}.
We insert $\mathbf{w}_1, \mathbf{w}_2, \mathbf{w}_3$ into the ODEs~\eqref{eq:TubODE1}--\eqref{eq:TubODE2} and append the auxiliary dynamic equations. Thus, the lifted discretized system is
\begin{align}
	\dot{\boldsymbol{\spe}} & = \A_{\boldsymbol{\spe}} \boldsymbol{\spe} + \mathbf{b}_\spe  u(t) - \Dam \ \boldsymbol{\spe} \odot  \mathbf{w}_1 ,\\
	\dot{\boldsymbol{\theta}} & = \A_{\boldsymbol{\theta}} \boldsymbol{\theta} + \mathbf{b}_\theta  u(t) + \mathcal{B}\Dam \ \boldsymbol{\spe} \odot  \mathbf{w}_1, \\
	\dot{\mathbf{w}}_1 & = \gamma \ \mathbf{w}_1 \odot \mathbf{w}_2 \odot \left [ A_2 \boldsymbol{\theta} + \mathbf{b}_\theta u(t) + \mathcal{B}\Dam \ \boldsymbol{\spe} \odot  \mathbf{w}_1 \right ],\\
	\dot{\mathbf{w}}_2 & = -2 \ \mathbf{w}_2 \odot \mathbf{w}_3 \odot \left [ A_2 \boldsymbol{\theta} + \mathbf{b}_\theta u(t) + \mathcal{B}\Dam \ \boldsymbol{\spe} \odot  \mathbf{w}_1 \right ], \\
	\dot{\mathbf{w}}_3 & = -\mathbf{w}_2 \odot \left [ A_2 \boldsymbol{\theta} + \mathbf{b}_\theta u(t) + \mathcal{B}\Dam \ \boldsymbol{\spe} \odot  \mathbf{w}_1 \right ].
\end{align}
The state of the lifted system is denoted as $\x = [ {\boldsymbol{\spe}}^\top, {\boldsymbol{\theta}}^\top, \mathbf{w}_1^\top, \mathbf{w}_2^\top, \mathbf{w}_3^\top ]^\top$.  We can write these equations as a quartic systems of ODEs as in Equation~\eqref{eq:quarticODE}  with $\B =  [\mathbf{b}_\spe^\top, \  \mathbf{b}_\theta^\top, \ \bzero ]^\top, \ \A = \text{diag}(\A_\psi, \A_\theta, \bzero_{3n})$ a block-diagonal matrix, and $\G^{(i)}, \N^{(i)} \in \mathbb{R}^{n\times n^i}$ being sparse matrices.
Given initial conditions $\boldsymbol{\theta}_0$ and $\boldsymbol{\spe}_0$, we find consistent initial conditions for $\mathbf{w}_1, \mathbf{w}_2, \mathbf{w}_3$ by using the definitions of the auxiliary variables in Equation~\eqref{eq:wDef}.	
This fourth-order polynomial ODE system is \textit{equivalent} to the the original ODE system \eqref{eq:TubODE1}--\eqref{eq:TubODE2} in that solutions $\boldsymbol{\spe}, \boldsymbol{\theta}$ for both systems are identical. However, the structure of the system is---as desired---polynomial, at the expense of increasing the discretization dimension from $2n$ to $5n$.

\subsection{Lifted Model 2: A Quadratic-Bilinear Differential-Algebraic System} \label{sec:TubQBDAE}
We further reduce the polynomial order of the system by lifting it to QB form. This requires introducing the following new dependent variables
\begin{equation} \label{eq:w4w5w6}
	\mathbf{w}_4 = \boldsymbol{\spe}\odot \mathbf{w}_1, \quad \mathbf{w}_5 = \mathbf{w}_2 \odot \mathbf{w}_3, \quad \mathbf{w}_6 = \mathbf{w}_1\odot \mathbf{w}_2.
\end{equation}
This time, we need algebraic constraints to represent the system, as further differentiation of the variables in~\eqref{eq:w4w5w6} would not result in a QB system, in turn requiring additional auxiliary variables.
With these new variables, the QB-DAE system is:
\begin{align}
	\dot{\boldsymbol{\spe}} & = \A_{\boldsymbol{\spe}} \boldsymbol{\spe} + \mathbf{b}_\spe  u(t) - \Dam \  \mathbf{w}_4, \label{eq:DAE1} \\
	\dot{\boldsymbol{\theta}} & = \A_{\boldsymbol{\theta}} \boldsymbol{\theta} + \mathbf{b}_\theta  u(t) + \mathcal{B}\Dam \  \mathbf{w}_4, \\
	\dot{\mathbf{w}}_1 & = \gamma \ \mathbf{w}_6\odot  \left [ A_2 \boldsymbol{\theta} + \mathbf{b}_\theta u(t) \right ] + \gamma \mathcal{B}\Dam \ \mathbf{w}_4 \odot  \mathbf{w}_6, \\
	\dot{\mathbf{w}}_2 & = -2 \ \mathbf{w}_5 \odot \left [ A_2 \boldsymbol{\theta} + \mathbf{b}_\theta u(t) \right ] -2  \mathcal{B} \Dam \ \mathbf{w}_4 \odot  \mathbf{w}_5,  \\
	\dot{\mathbf{w}}_3 & = -\mathbf{w}_2 \odot \left [ A_2 \boldsymbol{\theta} + \mathbf{b}_\theta u(t) \right ]  - \mathcal{B}\Dam \ \mathbf{w}_2 \odot  \mathbf{w}_4, \\
	0 & = \mathbf{w}_4 - \mathbf{w}_1 \odot \boldsymbol{\spe}, \label{eq:w4} \\
	0 & = \mathbf{w}_5 - \mathbf{w}_2 \odot \mathbf{w}_3 ,\label{eq:w5} \\
	0 & = \mathbf{w}_6 - \mathbf{w}_1 \odot \mathbf{w}_2 . \label{eq:w6}
\end{align}
The above system is a DAE of index 1. In other words, differentiating the algebraic constraints one time gives us an explicit ODE in terms of the other state variables.
We partition the state of the system into the dynamically evolving unconstrained states $\x_1$ and the states $\x_2$ that occur in the algebraic variables:
\begin{equation*}
	\x  = [\x _1^\top, \x _2^\top]^\top = [\underbrace{{\boldsymbol{\spe}}^\top, {\boldsymbol{\theta}}^\top, \mathbf{w}_1^\top, \mathbf{w}_2^\top, \mathbf{w}_3^\top}_{\x _1, \ \text{unconstrained}}, \ \underbrace{\mathbf{w}_4^\top, \mathbf{w}_5^\top, \mathbf{w}_6^\top }_{\x _2, \ \text{constrained}}]^\top.
\end{equation*}
The system~\eqref{eq:DAE1}--\eqref{eq:w6} can be written as a QB-DAE of the form~\eqref{eq:qb-dae} with matrices as in Equation~\eqref{eq:QBDAE-matrices} where $n_1 = 5n, n_2=3n$ and the mass matrix $\E_{11} = \I_{5n}$.  Moreover, the matrix $\B_1 =  [\mathbf{b}_\spe^\top, \  \mathbf{b}_\theta^\top, \ \bzero ]^\top $ and
$$
\A_{11} = \begin{bmatrix} \A_\psi & & \\ & \A_\theta & \\ & & \I_{3n}  \end{bmatrix}, \qquad
\A_{12} = \begin{bmatrix}  -\Dam \I_n & \bzero & \bzero \\ \mathcal{B}\Dam \I_n &\bzero & \bzero \end{bmatrix}.
$$
Here, $\A_{11}$ is the same as the matrix $\A$ in the quartic ODE of Section~\ref{sec:TubquarticODE}.

\subsection{Proper Orthogonal Decomposition for Quartic and Quadratic-Bilinear Differential-Algebraic Equations} \label{sec:TubPODdetails}

We compute ROMs of the quartic system and QB-DAE via projection onto POD basis vectors as described in Sections~\ref{sec:polySys} and \ref{sec:QB}.
We compute separate modes for each dependent variable. To illustrate this for the species concentration $\boldsymbol{\spe}$, let
\begin{equation}
	\Psi = [\boldsymbol{\spe}(t_0), \ \boldsymbol{\spe}(t_1), \ldots, \boldsymbol{\spe}(t_M)] \label{eq:snapshotMatrix}
\end{equation}
be the matrix of solution snapshots at equidistant times $t_i, \ i=1,\ldots, M$.
We compute the singular value decomposition of $\Psi = \mathbf{U} \mathbf{\Sigma} \mathbf{W}^\top$ and obtain the POD modes by taking the leading $r$ left singular vectors, $\V_\spe = \mathbf{U}(:, 1:r)$. Here, $r$ is chosen such that the system satisfies a certain accuracy level, as indicated by the decay in the singular values in $\mathbf{\Sigma}$.
The POD modes for the other dependent variables $\boldsymbol{\theta}$ and $\mathbf{w}_i, i =1,\ldots, 6$ are computed similarly, and stored in matrices $\V_\theta$ and  $\V_{\mathbf{w}_i}$ for $i =1,\ldots, 6$.

For the quartic system, the projection matrix is $\V = \text{blkdiag}(\V_\spe, \V_\theta, \V_{\mathbf{w}_1}, \V_{\mathbf{w}_2}, \V_{\mathbf{w}_3}) \in \mathbb{R}^{5n\times r}$ and used as in Equation~\eqref{eq:quarticROM} to obtain the quartic ROM.
For the QB-DAE system,  $\V_1 = \text{blkdiag}(\V_\spe, \V_\theta, \V_{\mathbf{w}_1}, \V_{\mathbf{w}_2}, \V_{\mathbf{w}_3}) \in \mathbb{R}^{5n\times r_1}$ and $\V_2= \text{blkdiag} (\V_{\mathbf{w}_4}, \V_{\mathbf{w}_5}, \V_{\mathbf{w}_6}) \in \mathbb{R}^{3n\times r_2}$ are the projection matrices used to obtain a QB-DAE ROM of the form~\eqref{eq:ROM_x1}--\eqref{eq:ROM_x2}.

\add{As illustrated in Section~\ref{sec:polySys}, the reduction of QB-ODEs or QB-DAEs does not require hyper-reduction. However, by using a projection matrix $\V_2\neq \mathbf{I}$ in Equation~\eqref{eq:ROM_x2}, we enforce the algebraic constraint---which encodes the part of the lifting transformation---only in the subspace $\V_2$. In that sense, the original nonlinearity is also approximated in our approach, but differently than in DEIM. 
	A similar statement holds for the QB-ODE case, where the auxiliary dynamics---again encoding the lifting transformation---are also projected onto $\V_2$ and thus introduce an approximation to the nonlinearity.}

\subsection{Numerical Results} \label{sec:numerics}
\begin{figure}[H]
	\begin{subfigure}{}
		\includegraphics[width=7cm]{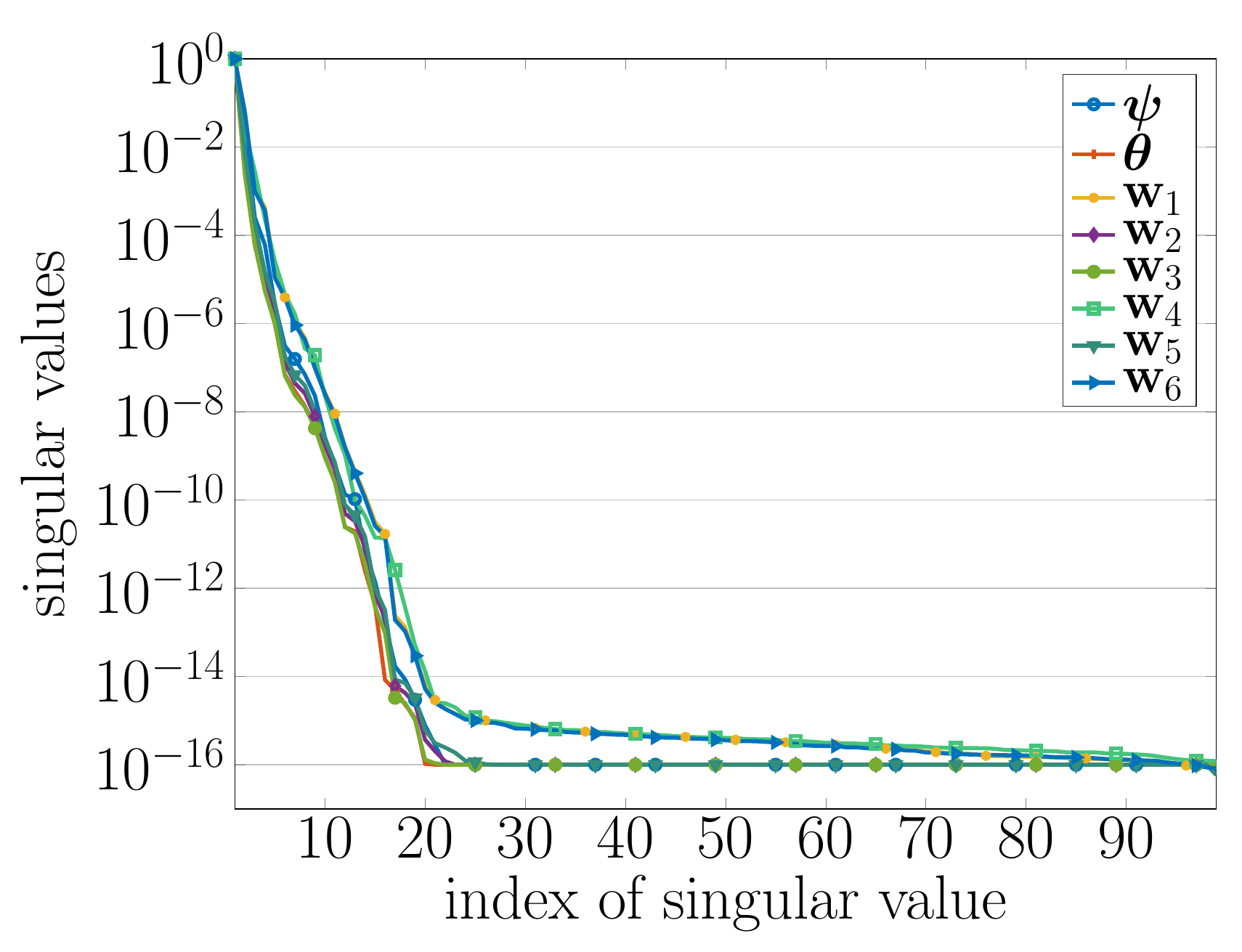}
	\end{subfigure}	
	\begin{subfigure}{}
		\includegraphics[width=7cm]{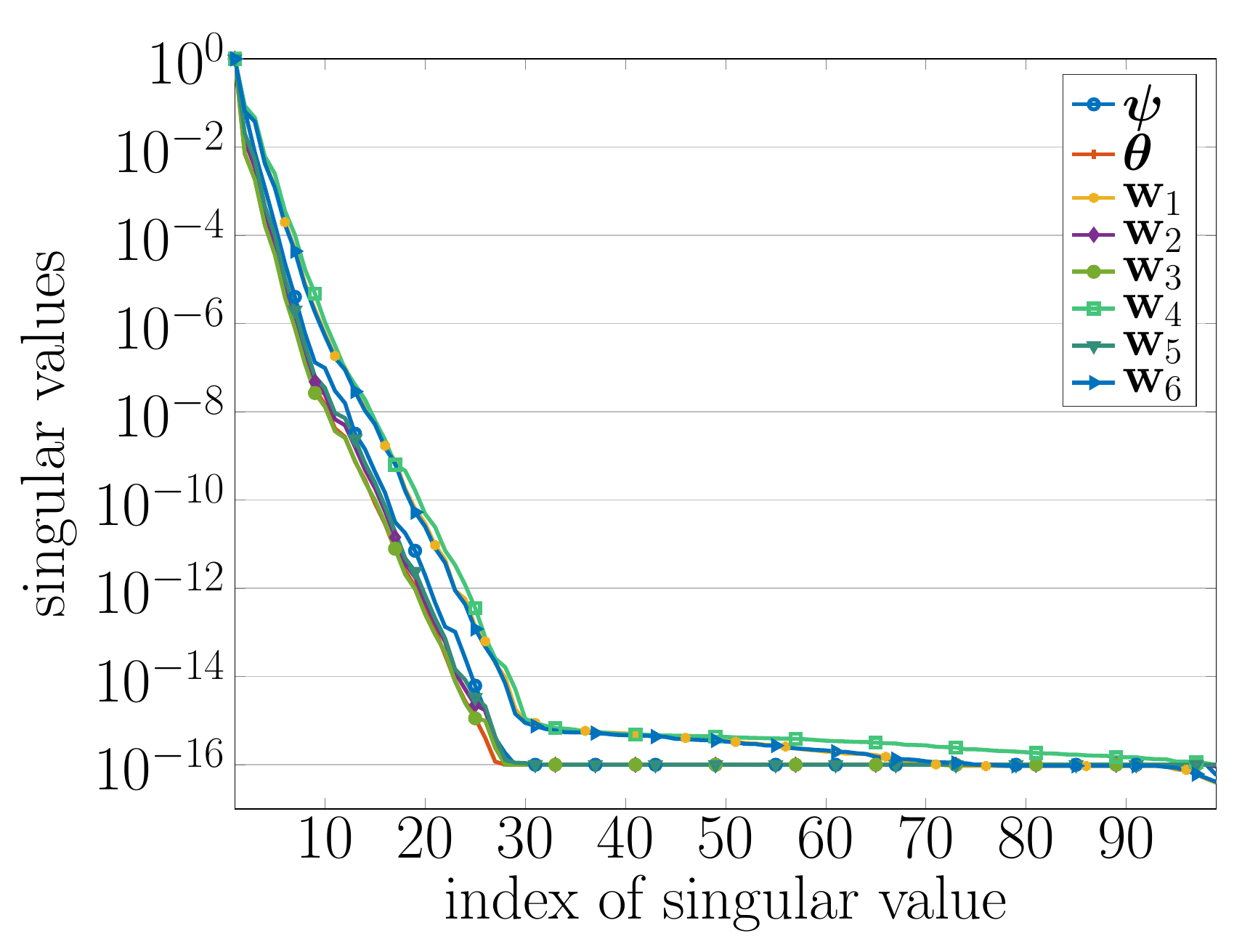}
	\end{subfigure}			
	\caption{Relative singular value decay of snapshot data for tubular reactor simulation.  Left: Damk{\"o}hler number $\Dam=0.162$ (stable case). Right: $\Dam=0.167$ (unstable case). The data according to $\mathbf{w}_1, \mathbf{w}_4, \mathbf{w}_6$ decays slower, as it relates to the exponential nonlinearity (see definition of those variables in equations~\eqref{eq:wDef} and \eqref{eq:w4w5w6}).}
	\label{fig:SVD_Da}
\end{figure}

\begin{figure}[H]
	\begin{subfigure}{}
		\includegraphics[width=7cm]{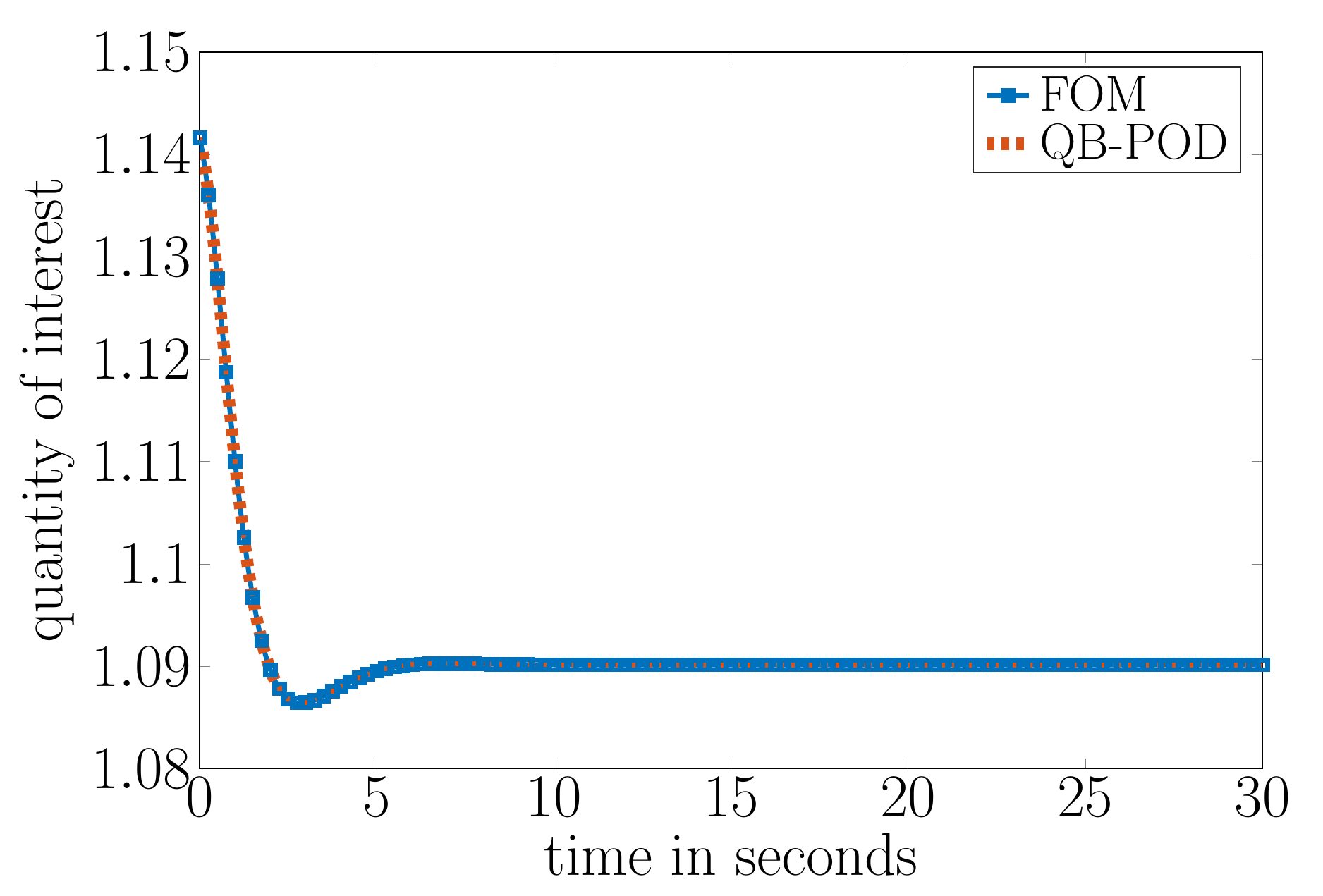}
	\end{subfigure}	
	\begin{subfigure}{}
		\includegraphics[width=7cm]{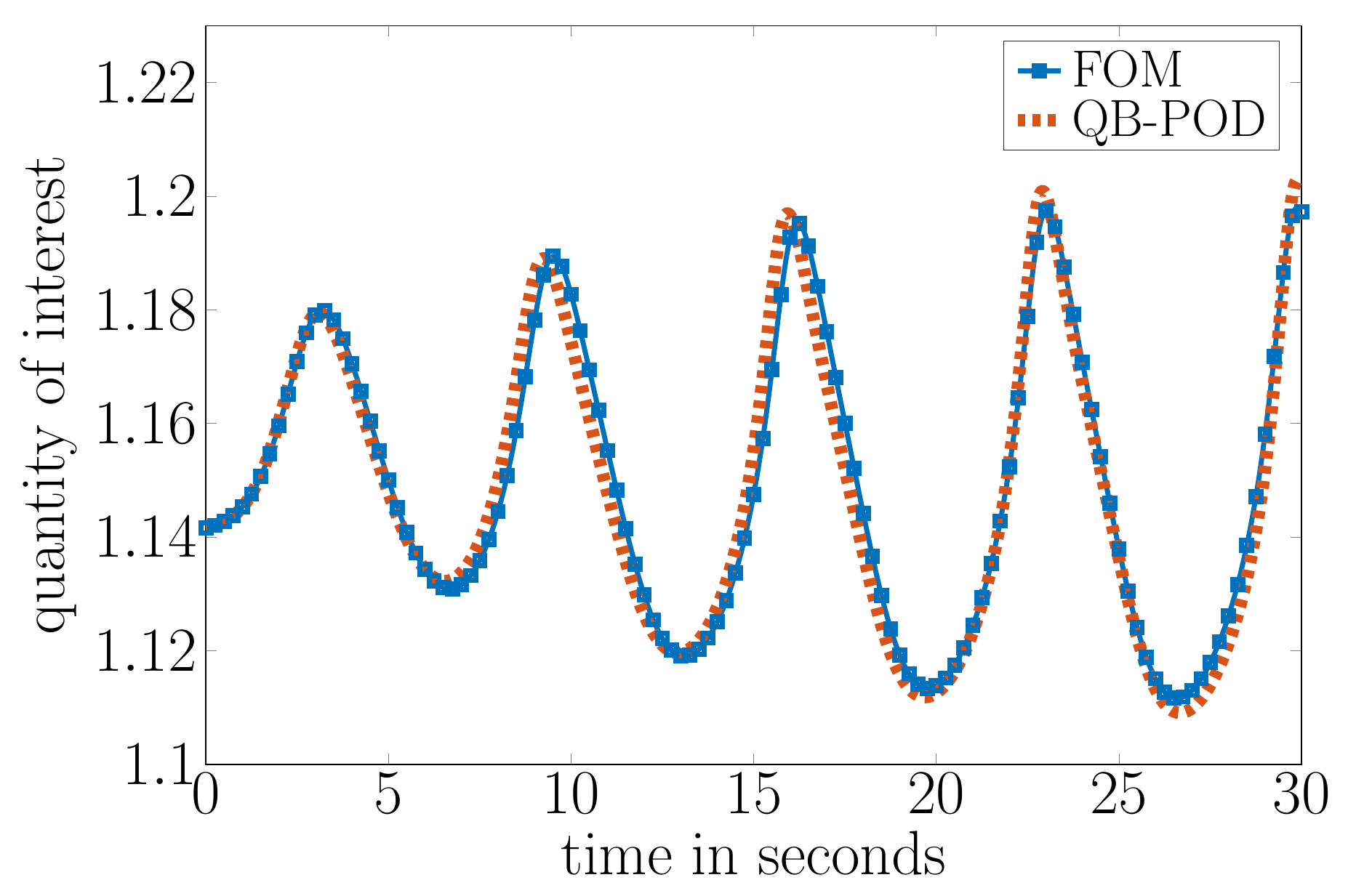}
	\end{subfigure}		
	\caption{Quantity of interest computed from FOM and QB-DAE ROM with $r_1 = 30, r_2=9$.  %
		Left: $\mathcal{D}=0.162$ (stable case), where the average relative state error from~\eqref{eq:state_error} is \add{$6.71\times 10^{-5}$}.
		Right: $\mathcal{D}=0.167$ (unstable case), where the average relative state error is \add{$8.95\times 10^{-3}$.}}
	\label{fig:TubPODoutput}
\end{figure}

We simulate the tubular reactor with parameters $Pe = 5, \gamma = 25, \mathcal{B} = 0.5, \beta = 2.5, \theta_{\text{ref}} \equiv 1$ until the final time of \add{$t_f=30s$} and record a snapshot every $\Delta t=0.01s$. The same initial conditions are used as in Ref.~\cite{zhou2012thesis}.
\add{For the computation of the POD basis, we only use snapshots until $t=20s$; thus, all the POD models in this section predict 50\% further past the training data}. 
Figure~\ref{fig:SVD_Da} shows the relative POD singular values for each variable. The left plot shows the stable case with $\Dam=0.162$ and the right plot the unstable regime with $\Dam=0.167$. The decay of the singular values is faster for all the variables in the case of low Damk{\"o}hler number, as is expected from the decaying stable dynamics.
Moreover, we see from Figure~\ref{fig:SVD_Da}, right, that the POD singular values for the variables $\boldsymbol{\theta}, \boldsymbol{\psi}, \mathbf{w}_2, \mathbf{w}_3$, and  $\mathbf{w}_5$ all decay similarly. A slower decay of the POD singular values corresponding to the variables $\mathbf{w}_1, \mathbf{w}_4$, and  $\mathbf{w}_6$ is observed compared to the decay of $\boldsymbol{\theta}, \boldsymbol{\psi}, \mathbf{w}_2, \mathbf{w}_3$, and $\mathbf{w}_5$. Note that $\mathbf{w}_4, \mathbf{w}_6$ are related to $\mathbf{w}_1$, see Equation~\eqref{eq:w4w5w6}, hence their similar decay in POD singular values.
The quantity of interest for this example is the temperature oscillation at the reactor exit, $\theta(s=1,t)$. Figure~\ref{fig:TubPODoutput} shows the quantity of interest predicted using the FOM and the QB-DAE ROM, which was generated via POD as described above with $r_1=30$ and $r_2=9$ basis functions. The case of stable dynamics is plotted on the left, and dynamics with limit-cycle oscillations are shown on the right plot. The QB-DAE ROMs are accurate in both cases and reproduce the limit-cycle amplitude and oscillations well. We note that, as is often the case with POD ROMs, not all choices of basis size yield satisfactory ROMs. \add{For instance, in this strongly nonlinear example, we found that for a fixed number of modes, certain selections of modes used in $\V_\spe, \V_\theta, \V_{\mathbf{w}_1}, \V_{\mathbf{w}_2}, \V_{\mathbf{w}_3}$ gave better results than others.} 

\begin{figure}[H]
	\begin{subfigure}{}
		\includegraphics[width=7cm]{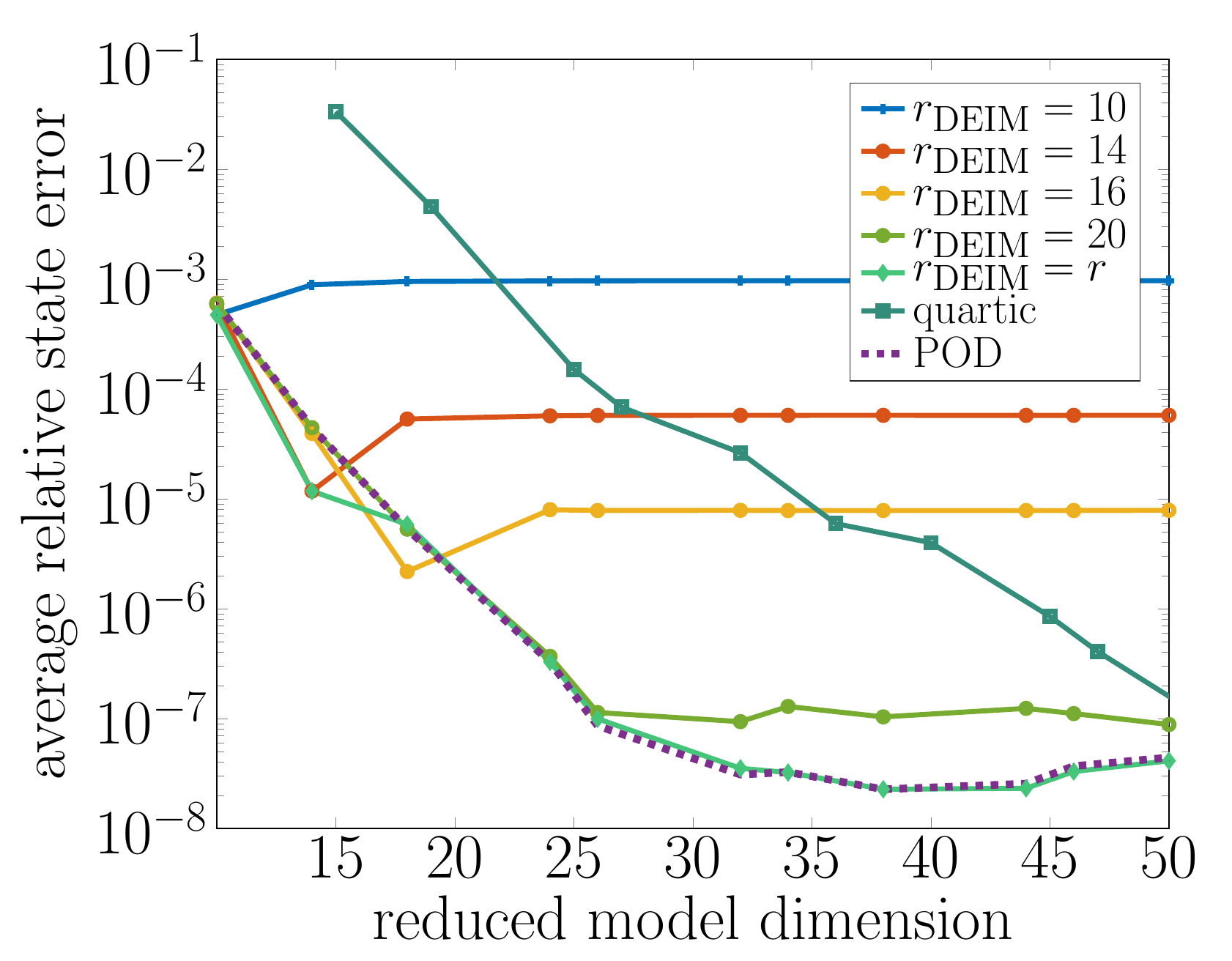}
	\end{subfigure}		
	\begin{subfigure}{}
		\includegraphics[width=7cm]{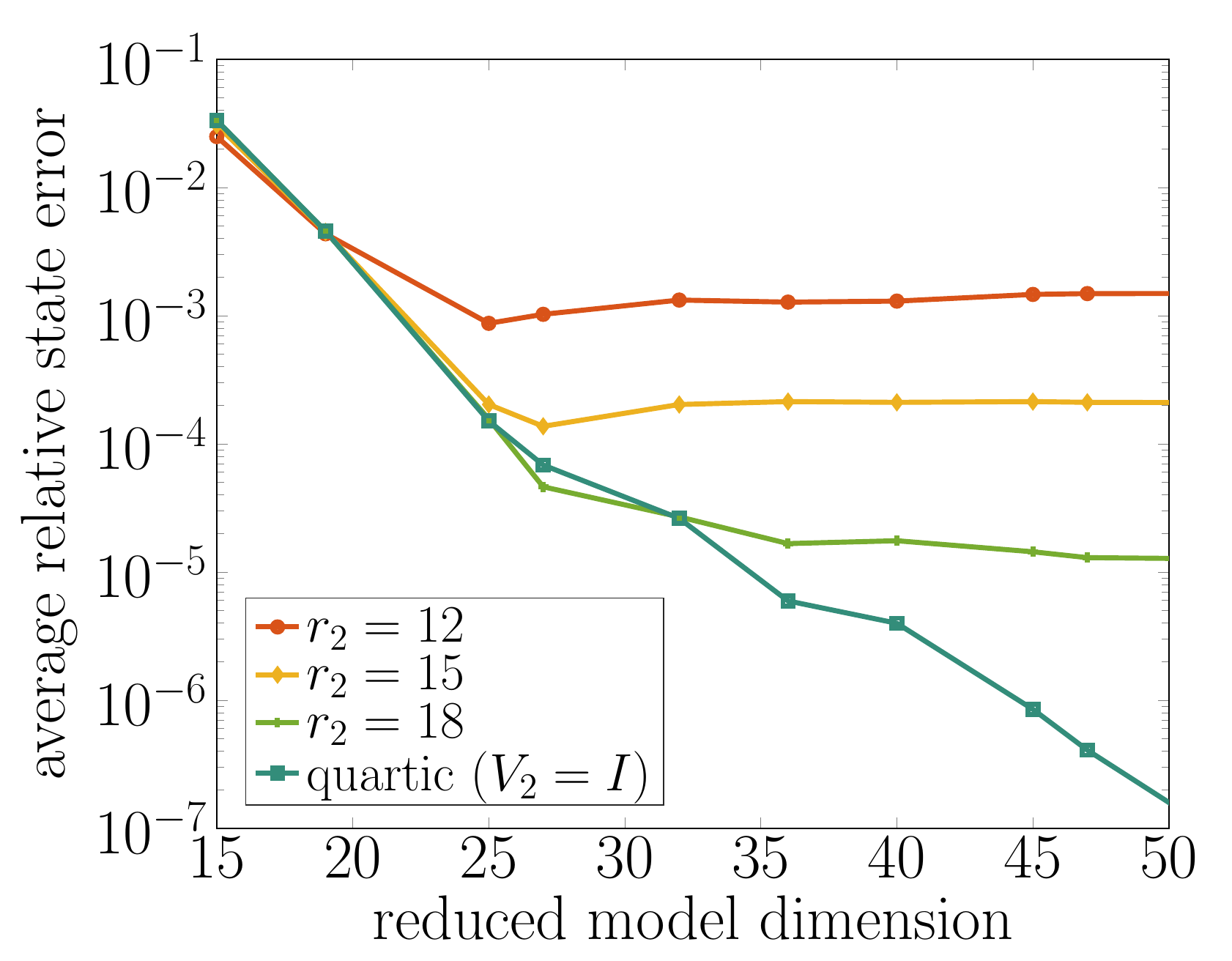}
	\end{subfigure}	
	\caption{Comparison of average relative state errors as ROM dimension increases for $\mathcal{D}=0.167$. Left:  We compare DEIM approximations with \add{a fixed number of DEIM interpolation points, a POD DEIM model with $r_\text{DEIM} = r$,} the quartic POD ROM, and a POD ROM of the original model. Right: Comparison of the QB-DAE ROM approximation for different approximations of the constrained states $\x_2 \approx \V_2 \widehat{\x}_2$, where $r_2$ is the number of basis functions in $\V_2$.  }
	\label{fig:TubErrComp}
\end{figure}

We compute average relative state errors from \add{$n_t=3000$} state snapshots of FOM and ROM solutions as
\begin{equation} \label{eq:state_error}
	\epsilon(r_1,r_2) = \frac{1}{n_t} \sum_{i=1}^{n_t} \Vert \x(t_i)  - \x^\text{ROM}(t_i) \Vert /  \Vert \x(t_i) \Vert,
\end{equation}
where $\x(t_i) = [ \boldsymbol{\spe}^\top,  \boldsymbol{\theta}^\top]^\top (t_i)$ is the solution of the FOM at time step $t_i$ and likewise $\x^\text{ROM}(t_i)$ is the ROM solution of the original variables at step $t_i$, i.e., we only compare the approximation in the original state variables $\boldsymbol{\spe}$ and $ \boldsymbol{\theta}$. The error is given as a function of $r_1$ and $r_2$, which are the numbers of POD modes used in $\V_1$ and $\V_2$ in Equation~\eqref{eq:defV}.
Figure~\ref{fig:TubErrComp}, left, shows the error $\epsilon(r_1,r_2)$  for different ROMs: Four POD-DEIM reduced models of the FOM~\eqref{eq:TubODE1}--\eqref{eq:TubODE2} with $r_{\text{DEIM}}=10, 14, 16, 20$ DEIM interpolation points; \add{a POD-DEIM model that increases the DEIM interpolation points with the reduced-dimension, i.e., we have $r_{\text{DEIM}}=r$}; a standard POD reduced model; and the quartic ROM from Equation~\eqref{eq:quarticROM}.
The POD approximation provides the lower bound on the error, as it directly evaluates the full nonlinear right-hand side of~\eqref{eq:TubODE1}--\eqref{eq:TubODE2}, which scales in computational complexity with the full state dimension $2n$. Even though the POD model is accurate, it is not computationally feasible, and is shown only for reference.
The DEIM approximations are less accurate than the POD model, but increase in accuracy when more DEIM interpolation points are used. As is typically observed with POD-DEIM reduced models (e.g., Ref.~\cite{deim2010}), the interpolation error dominates after some time \add{for a fixed number of interpolation points}, and so the model cannot improve further as more basis vectors are added. \add{The POD-DEIM model with $r_{\text{DEIM}}=r$ however approximates the original POD model well.}
The quartic ROM \add{likewise} does not suffer from the limitation of hyper-reduction interpolation error and increases in accuracy as further basis functions $r_1$ are added.
Figure~\ref{fig:TubErrComp}, right, shows the influence of the approximation of the constrained states on the accuracy of the QB-DAE reduced model. 
As mentioned in Section~\ref{sec:QB}, the case $\V_2 = \I$ leads to the quartic ODE. We compare this with the three cases of $r_2= 12, 15, 18$. Figure~\ref{fig:TubErrComp}, right, shows the state errors plotted against $r_1$, the approximation dimension of the dynamic variables $\x_1$. We observe a similar trend by increasing the approximation $\x_2 \approx \V_2 \widehat{\x}_2$ as compared to increasing the DEIM interpolation points. The better the approximation of the constrained states $\x_2$, the more accurate the corresponding QB-DAE ROM. 

\section{Conclusions} \label{sec:conclusion}

The approach to first lift a nonlinear dynamical system via auxiliary variables and then reduce the structured problem presents an alternative to the state-of-the-art in nonlinear reduced-order modeling. The nonlinear partial differential equations arising in many aerospace systems of interest can be lifted to have polynomial form; lifting to a system of quadratic-bilinear DAEs is shown to have particular promise.
We derived multi-step lifting transformations for a strongly nonlinear tubular reactor model.
The numerical results show that the lifting approach together with structure-preserving POD-based model reduction is competitive with the state-of-the-art POD-DEIM nonlinear model reduction approach for the chosen examples. The lifting approach has the added advantage of introducing no additional approximation in the reduced model nonlinear terms; this comes at the cost of the extra up-front work to derive the lifted system, including the matrices and tensors that represent the lifted dynamics. Another advantage is that the structured polynomial form of the lifted systems
holds significant promise for building more rigorous analysis of ROM stability and error behavior, especially for quadratic-bilinear systems.
The results also highlight a potential drawback in that the introduction of auxiliary variables increases the dimension of the state and also tends to increase the number of POD basis vectors needed to achieve an acceptable error. This drawback could potentially be addressed by using nonlinear projection subspaces in place of the linear POD subspaces, which is particularly viable since the lifting transformations are known.

\aiaaappendix{Details on Solving the QB-DAE ROM} \label{sec:appendix}
This section expands on the details regarding the solution of the QB-DAE ROM in equations~\eqref{eq:ROM_x1}--\eqref{eq:ROM_x2}.
One can either solve those equations with specific DAE solvers (such as \texttt{ode15s} in Matlab), or we can simulate the DAE by inserting $\widehat{\x}_2$ into the dynamic equations. When doing so, we can speed up the simulations by efficiently pre-computing the matricized tensor as follows:
\begin{align*}
	\widehat{\HH}_1 \left (
	\begin{bmatrix}
		\widehat{\x}_1 \\ \widehat{\x}_2 \end{bmatrix} \otimes \begin{bmatrix}
		\widehat{\x}_1 \\ \widehat{\x}_2 \end{bmatrix} \right )
	& =
	\widehat{\HH}_1\left (
	\begin{bmatrix} \I_{r_1} &  \\ & \widehat{\HH}_2 \end{bmatrix}
	\begin{bmatrix} \widehat{\x}_1 \\ \widehat{\x}_1 \otimes \widehat{\x}_1 \end{bmatrix}
	\right )
	\otimes
	\left (
	\begin{bmatrix} \I_{r_1} &  \\ & \widehat{\HH}_2 \end{bmatrix}
	\begin{bmatrix} \widehat{\x}_1 \\ \widehat{\x}_1 \otimes \widehat{\x}_1 \end{bmatrix}
	\right ) \\
	& =
	\widehat{\HH}_1 \left (
	\begin{bmatrix} \I_{r_1} &  \\ & \widehat{\HH}_2 \end{bmatrix}
	\otimes
	\begin{bmatrix} \I_{r_1} &  \\ & \widehat{\HH}_2 \end{bmatrix}
	\right )
	\left (
	\begin{bmatrix} \widehat{\x}_1 \\ \widehat{\x}_1 \otimes \widehat{\x}_1 \end{bmatrix}
	\otimes
	\begin{bmatrix} \widehat{\x}_1 \\ \widehat{\x}_1 \otimes \widehat{\x}_1 \end{bmatrix}
	\right ) \\
	& =: \widetilde{\HH}_1 \begin{bmatrix} \widehat{\x}_1 \otimes \widehat{\x}_1 \\ \widehat{\x}_1 \otimes \widehat{\x}_1 \otimes \widehat{\x}_1 \\ \widehat{\x}_1 \otimes \widehat{\x}_1 \otimes \widehat{\x}_1 \otimes \widehat{\x}_1 \end{bmatrix}
\end{align*}
with $\widetilde{\HH}_1 \in \mathbb{R}^{r_1 \times (r_1^2 + r_1^3 + r_1^4)}$ and where the second equality follows from properties of the Kronecker product, i.e, $\A \C \otimes \B \mathbf{D} = (\A \otimes \B)(\C \otimes \mathbf{D})$.
Thus, we obtain the ODE
\begin{align*}
	\dot{\widehat{\x} }_1 & = \widehat{\A}_{11} \widehat{\x} _1 + \widehat{\A}_{12} \widehat{\HH}_2 (\widehat{\x}_1\otimes \widehat{\x}_1) + \widehat{\B}_1 u + \widetilde{\HH}_1 \begin{bmatrix} \widehat{\x}_1 \otimes \widehat{\x}_1 \\ \widehat{\x}_1 \otimes \widehat{\x}_1 \otimes \widehat{\x}_1 \\ \widehat{\x}_1 \otimes \widehat{\x}_1 \otimes \widehat{\x}_1 \otimes \widehat{\x}_1 \end{bmatrix} + \widehat{\N}_{11} \widehat{\x}_1 u + \widehat{\N}_{12} \widehat{\HH}_2 (\widehat{\x}_1\otimes \widehat{\x}_1) u  .
\end{align*}	
The matrix products $\widehat{\A}_{12} \widehat{\HH}_2$ and $\widehat{\N}_{12} \widehat{\HH}_2$	can be pre-computed offline for faster online computation. 	

\section*{Funding Sources}
This work was supported in part by the Air Force Center of Excellence on Multi-Fidelity Modeling of Rocket Combustor Dynamics under award FA9550-17-1-0195.

\section*{Acknowledgments}
The authors thank Dr. Pawan Goyal for sharing code related to the FitzHugh-Nagumo problem.

\bibliographystyle{abbrv}
\bibliography{KramerWillcox_nonlinearModelReductionLiftingPOD.bib}

\end{document}
